\documentclass[10pt]{article}

\usepackage{amsmath}
\usepackage{mathtools}
\usepackage{amssymb}
\usepackage{latexsym}
\usepackage{enumitem}
\usepackage{mathrsfs}
\usepackage{comment}
\usepackage{pifont}
\usepackage{enumitem}
\usepackage{bbold}
\usepackage{graphicx}
\usepackage{color}
\usepackage{hyperref}
\hypersetup{
    colorlinks=true,
    linkcolor=blue,
    filecolor=magenta,      
    urlcolor=cyan,
    citecolor=magenta
}

\newtheorem{assumption}{Assumption}

\def\qed{ \ \vrule width.2cm height.2cm depth0cm\smallskip}
\newenvironment{proof}{\noindent {\bf Proof.\/}}{$\qed$\vskip 0.1in}

\newcommand{\ol}{\overline}
\newcommand{\ul}{\underline}
\newcommand{\eps}{\varepsilon}

\newcommand{\ba}{\begin{array}}
\newcommand{\ea}{\end{array}}
\newcommand{\be}{\begin{equation}}
\newcommand{\ee}{\end{equation}}
\newcommand{\bea}{\begin{eqnarray}}
\newcommand{\eea}{\end{eqnarray}}
\newcommand{\beaa}{\begin{eqnarray*}}
\newcommand{\eeaa}{\end{eqnarray*}}

\def\dbD{\mathbb{D}}
\def\dbE{\mathbb{E}}
\def\dbF{\mathbb{F}}

\def\dbI{\mathbb{I}}

\def\dbL{\mathbb{L}}

\def\dbP{\mathbb{P}}
\def\dbR{\mathbb{R}}

\def\dbT{\mathbb{T}}

\def\dbZ{\mathbb{Z}}

\def\a{\alpha}
\def\b{\beta}
\def\g{\gamma}
\def\d{\delta}
\def\e{\varepsilon}
\def\z{\zeta}

\def\l{\lambda}

\def\si{\sigma}
\def\t{\tau}
\def\f{\varphi}

\def\D{\Delta}

\def\O{\Omega}
\def\cA{{\cal A}}
\def\cB{{\cal B}}

\def\cF{{\cal F}}

\def\cK{{\cal K}}
\def\cL{{\cal L}}
\def\cM{{\cal M}}
\def\cN{{\cal N}}

\def\cP{{\cal P}}

\def\cW{{\cal W}}

\def\no{\noindent}

\def\ms{\medskip}
\def\bs{\bigskip}
\def\q{\quad}
\def\qq{\qquad}

\def\pa{\partial}
\def\cd{\cdot}

\def\tr{\hbox{\rm tr}}

\def\qed{ \hfill \vrule width.25cm height.25cm depth0cm\smallskip}

\newcommand{\basa}{\begin{assumption}}
\newcommand{\easa}{\end{assumption}}

\newcommand{\bas}{\begin{assum}}
\newcommand{\eas}{\end{assum}}

\def\pa{\partial}

 \def\cd{\cdot}

\def\supp{\hbox{\rm Supp$\,$}}

\def\tr{\hbox{\rm tr$\,$}}

\def\dis{\displaystyle}

\def\1{{\bf 1}}

\def\:{\!:\!}
\def\reff#1{{\rm(\ref{#1})}}
\def \proof{{\noindent \bf Proof\quad}}

\def \bS{{\bf S}}
\def \bQ{{\bf Q}}
\def \bm{{\bf m}}

at 9pt

\begin{document}

\newtheorem{thm}{Theorem}[section]
\newtheorem{lem}[thm]{Lemma}
\newtheorem{cor}[thm]{Corollary}
\newtheorem{prop}[thm]{Proposition}
\newtheorem{rem}[thm]{Remark}
\newtheorem{eg}[thm]{Example}
\newtheorem{defn}[thm]{Definition}
\newtheorem{assum}[thm]{Assumption}

\numberwithin{equation}{section}

\title{\bf{Viscosity solutions for obstacle problems on Wasserstein space\footnote{The first two authors are grateful for the financial support from the Chaires FiME-FDD and Financial Risks of the Louis Bachelier Institute. The third author is supported in part by NSF grant DMS-1908665.}}}
\author{Mehdi Talbi\footnote{CMAP, \'{E}cole polytechnique, France, mehdi.talbi@polytechnique.edu} \quad Nizar Touzi\footnote{CMAP, \'{E}cole polytechnique, France, nizar.touzi@polytechnique.edu}  \quad Jianfeng Zhang\footnote{Department of Mathematics, University of Southern California, United States, jianfenz@usc.edu. }}
\date{\today}

\maketitle

\begin{abstract}
This paper is a continuation of our accompanying paper \cite{TTZ}, where we characterized the mean field optimal stopping problem by an obstacle equation on the Wasserstein space of probability measures, provided that the value function is smooth. Our purpose here is to establish this characterization under weaker regularity requirements. We shall define a notion of viscosity solutions for such equation, and prove existence, stability, and comparison principle. 
\end{abstract}

\no{\bf MSC2020.} 60G40, 35Q89,  49N80, 49L25, 60H30.

\vspace{3mm}
\no{\bf Keywords.} Mean field optimal stopping, obstacle problems, viscosity solutions

\section{Introduction}

In our previous paper \cite{TTZ}, we characterized the so-called mean field optimal stopping problem by a dynamic programming equation on the Wasserstein space, that we call obstacle equation on Wasserstein space by analogy with the equation corresponding to the standard optimal stopping problem (see e.g. El Karoui \cite{EK} or Shiryaev \cite{Shiryaev}). More precisely, we proved that the value function of our optimization problem is the unique solution of the obstacle equation on the Wasserstein space, provided it has $C^{1,2}$ regularity (in an appropriate sense). 
We note that, besides its obvious connection with multiple stopping problems over a large interacting particle system, this obstacle equation provides a convenient tool for many time inconsistent optimal stopping problems. We also remark that, our mean field optimal stopping problem has quite different structure than the mean field games of optimal stopping.  

However, as in the case of the standard optimal stopping problems, one can rarely expect a classical solution for the obstacle equations. In particular, the infinite dimensionality of the space of measures makes the regularity requirement even harder to meet. Our goal of this paper is thus to develop a viscosity solution theory for the obstacle problem on the Wasserstein space, which as well-known requires much weaker regularities. 

There have been some serious efforts on viscosity solutions of nonlinear partial differential equations on the Wasserstein space. We first mention the paper by Cardaliaguet \& Quincampoix \cite{CarQui}, which considered a first order Hamilton-Jacobi-Isaacs equation on Wasserstein space arising from deterministic zero-sum games with random initial conditions. The comparison principle for viscosity solutions was established by combining the doubling variables argument with Ekeland's variational principle. We may also mention the work of Gangbo, Nguyen \& Tudorascu \cite{GNT} and Jimenez, Marigonda \& Quincampoix \cite{JMQ}, who also define a notion of viscosity solutions for Hamilton-Jacobi equations by using sub-differentials.
Another approach followed by several authors (see e.g.\ Pham \& Wei \cite{PW}) consists in exploiting Lions' idea \cite{Lions} of lifting the functions on the Wasserstein space into functions on the Hilbert space of random variables and then using the existing viscosity theory on Hilbert spaces (see e.g.\ Lions \cite{Lions1, Lions2, Lions3} and Fabbri, Gozzi \& Swiech \cite{FGS}). More recently, Cosso, Gozzi, Kharroubi, Pham \& Rosestolato \cite{CGKPR} defined viscosity solutions for Hamilton-Jacobi equations by requiring the global extrema on the Wasserstein space for the tangency property of the test functions.

In the context of mean field control problems in a path dependent setting, Wu \& Zhang \cite{WZ} proposed a notion of viscosity solutions for parabolic equations on the Wasserstein space, inspired from Ekren, Keller, Ren, Touzi \& Zhang \cite{EKTZ,ETZ1,ETZ2, RTZ2}. Note that the natural idea which consists in taking Wasserstein balls for the viscosity neighborhood (as in Carmona \& Delarue \cite{CD2}) leads to serious difficulties as the Wasserstein ball is in general not compact. Instead, \cite{WZ} restricted the viscosity neighborhood of some point $(t,\mu)$ (where $t$ is a time and $\mu$ a measure) to the compact set of all possible laws of the controlled state process starting from this point. Another remarkable work by Burzoni, Ignazio, Reppen \& Soner \cite{BIRS}, in the context of mean field control of jump-diffusions, restricted the viscosity neighborhood in another way, so as to guarantee compactness. They proved a comparison result by the doubling variables argument. To do this, they succeeded in constructing a smooth metric which serves as a test function, but unfortunately restricts the scope of the method to the case when the coefficients of the controlled dynamics do not depend on the space variable. 

We shall follow the approach of \cite{WZ}. We consider the joint law of $(X_{\t \wedge t}, \1_{\{\t\ge t\}})$ as the variable of the value function, where $X$ is the state process and $\t$ is the stopping time.  As in \cite{WZ} we define viscosity solutions by using the set of such laws over all stopping times $\t$.  This neighborhood set of laws, for a given initial condition, is compact under Wasserstein distance and thus is desirable for the viscosity theory. We show that, under natural conditions, the value function of the mean field optimal stoping problem is indeed the unique viscosity solution of the corresponding obstacle equation on Wasserstein space. We shall also establish the stability and the comparison principle for the viscosity solutions. To prove the latter, one key ingredient is a smooth mollifier for continuous functions on the Wasserstein space, introduced by Mou \& Zhang \cite{MZ}.  However, to obtain some uniform estimates of the smooth mollifier which are needed in our proof of comparison principle,  as in \cite{MZ} we require the coefficients to be Lipschitz continuous under the $1$-Wasserstein distance, rather than the more natural  $2$-Wasserstein distance.

As applications of our viscosity theory, we invest several time inconsistent optimal stopping problems, including problems related to mean variance, probability distortion, and expected shortfall. By considering the law (instead of the value) of the stopped state process as the variable, we show that the value functions are indeed the unique viscosity solution to the corresponding obstacle equation on the Wasserstein space. Moreover, our results can be easily extended to the infinite horizon case.   

The paper is organized as follows. In \S\ref{sect-obstacle}, we present the mean field optimal stopping problem, the corresponding dynamic programming equation and some of its elementary properties. \S\ref{sect-viscosity} is the main section, where we propose our definition of viscosity solutions and prove the main results. \S\ref{sect-examples} is devoted to several applications. 
Finally, we prove some technical results in the appendix.

\vspace{3mm}
\no {\bf Notations.} We denote by $\cP(\O,\cF)$ the set of probability measures on a measurable  space $(\O,\cF)$, and $\cP_2(\O,\cF) := \{ m \in \cP(\O,\cF) : \int_\O d(x_0,x)^2 m(dx) < \infty \}$ for some $x_0 \in \O$, where $d$ is a metric on $\O$. $\cP_2(\O,\cF)$ is equipped with the corresponding $2$-Wasserstein distance $\cW_2$. When $(\O,\cF) = (\dbR^d,\cB(\dbR^d))$, we simply denote them as $\cP(\dbR^d)$ and $\cP_2(\dbR^d)$. For a random variable $Z$ and a probability $\dbP$, we denote by $\dbP_Z:=\dbP\circ Z^{-1}$ the law of $Z$ under $\dbP$. 
For vectors $x, y\in \dbR^n$ and matrices $A, B\in \dbR^{n\times m}$, denote $ x\cd y:=\sum_{i=1}^n x_iy_i$  and $A:B:= \tr(A B^\top)$. We shall also write “USC" (resp. “LSC") “upper (resp. lower) semi-continuous".

\section{The obstacle problem on Wasserstein space}\label{sect-obstacle}
\subsection{Formulation}

Let $T < \infty$ be fixed, and $\O := C^0([-1,T],\dbR^d) \times \dbI^0([-1,T])$ the canonical space, where: 

$\bullet$ $C^0([-1,T],\dbR^d)$ is the set of continuous paths from $[-1,T]$ to $\dbR^d$, constant on $[-1,0)$; 

$\bullet$ $\dbI^0([-1,T])$ is the set of non-increasing and càdlàg maps from $[-1,T]$ to $\{0,1\}$,  constant 

~ on $[-1,0)$, and ending with value $0$ at $T$.
 
 \no We equip $\O$ with the Skorokhod distance, under which it is a Polish space. Note that the choice of the extension to $-1$ is arbitrary,  the extension of time  to the left of the origin is only needed to allow for an immediate stop at time $t=0$. 
 
We denote $Y := (X,I)$ the canonical process, with state space $\mathbf{S}:=\dbR^d\times\{0,1\}$, its canonical filtration $\dbF=\dbF^Y= (\cF_t)_{t \in [-1,T]}$, and the corresponding jump time of the survival process $I$:
\bea
\label{tau}
\t := \inf\{t \ge 0 : I_t = 0\},  \ \mbox{so that $I_t := I_{0-}\1_{t < \t}$ for all $t \in [-1,T]$.} \nonumber
\eea 
By the càdlàg property of $I$, $\t$ is an $\dbF-$stopping time. 

 Let $(b,\si,f): [0,T] \times \dbR^d \times \cP_2(\mathbf{S}) \rightarrow \dbR^d \times \dbR^{d\times d} \times \dbR$ with $\si$ taking values in non-negative matrices, and $ g : \cP_2(\dbR^d) \rightarrow \dbR$.  In the following assumption, which will always be in force throughout the paper, $\cP_2(\mathbf{S})$ is equipped with the $\cW_2$ distance. 
\begin{assum}
\label{assum-bsig}
\no{\rm (i)} $b, \si$ are continuous in $t$, and uniformly Lipschitz continuous in $(x, m)$.

\no{\rm (ii)} $f$ is Borel measurable and has quadratic growth in $x\in \dbR^d$, and 
\bea
\label{F}
F(t,m) := \int_{\dbR^d} f(t,x,m)m(dx,1)
~\mbox{is continuous on}~[0, T]\times \cP_2(\bS).
\eea
\no{\rm (iii)}  $g$ is upper-semicontinuous and locally bounded; and extended to $\cP_2(\bS)$ by $g(m) := g(m(\cd, \{0,1\}))$.
\end{assum}
 
Introduce the dynamic value function
 \begin{equation}\label{weakoptstop}
 V(t,m) := \underset{\dbP \in \cP(t,m)}{\sup} \Big\{ \int_t^T  F(r,  \dbP_{Y_r})dr  + g(\dbP_{Y_T}) \Big\}, \quad \mbox{$(t,m) \in [0,T] \times \cP_2(\mathbf{S})$}.
 \end{equation}
 Here $\cP(t,m)$ is the set of probability measures $\dbP$ on $(\O,\cF_T)$ s.t. $\dbP_{Y_{t-}} = m$,  the paths $s\in [-1, t) \to Y_s$ are constants, $\dbP$-a.s., and  the processes:
 \bea\label{martingalepb}
M_. := X_. - \int_t^. b(r,X_r,\dbP_{Y_r})I_rdr \q \mbox{and} \q M_. M_.^\intercal -  \int_t^. \si^2(r,X_r,\dbP_{Y_r})I_rdr
\eea
are $\dbP-$martingales on $[t,T]$, that is, for some $\dbP-$Brownian motion $W^\dbP$, 
\begin{equation}\label{asympt}
X_s = X_t + \int_t^s b(r, X_r, \dbP_{Y_r})I_r dr + \int_t^s \sigma(r, X_r, \dbP_{Y_r})I_r dW_r^\dbP , \ I_s = I_{t-} \1_{s < \t}, \ \dbP-\mbox{a.s.} \nonumber
\end{equation}
A special element of $\cP(t, m)$ is $\bar \dbP=\bar \dbP^{t,m}$ under which $X$ is unstopped. That is,
\bea
\label{barP}
X_s = X_t + \int_t^s b(r, X_r, \bar \dbP_{Y_r}) I_r dr + \int_t^s \si(r, X_r, \bar \dbP_{Y_r})I_r  dW^{\bar \dbP}_r,~ I_s = I_{t-}\1_{[t, T)}(s),~ \bar \dbP \mbox{-a.s.} 
\eea
Note that $Y_. = Y_{. \wedge \t}$, and in particular $Y_T = Y_\t,$ $\dbP-$a.s. Moreover, from the definition of $F$ in \eqref{F}, we have $\int_t^T F(r,\dbP_{Y_r})dr=\dbE^\dbP\int_t^\tau f(r,X_r,\dbP_{Y_r})dr$.

We recall from our first paper \cite{TTZ} that $\cP(t,m)$ is compact under the Wasserstein distance $\cW_2$, and thus existence holds for the mean field optimal stopping problem \eqref{weakoptstop}.  Furthermore, we have the dynamic programming principle (DPP for short): for any $s\in [t, T]$,
\bea\label{weakDPP}
 V(t,m) = \!\!\sup_{\dbP \in \cP(t,m)} \Big\{ \int_t^s F(r,  \dbP_{Y_r})dr   + V(s, {\dbP}_{Y_{s-}}) \Big\}  = \!\!\sup_{\dbP \in \cP(t,m)} \Big\{ \int_t^s  F(r, \dbP_{Y_r})dr   + V(s, {\dbP}_{Y_{s}}) \Big\}.
\eea

\subsection{Differential calculus}

We next recall some differential calculus tools on the Wasserstein space. We say that a function $u:\cP_2(\bS) \to \dbR$ has a functional linear derivative $\delta_m u:\cP_2(\bS)\times\bS \to \dbR$ if
 \beaa
u(m')-u(m) 
= 
\int_0^1 \int_{\bS} \d_m u(\l m' + (1-\l)m, y)(m'-m)(dy)d\l
&\mbox{for all}&
m,m'\in\cP_2(\bS),
\eeaa
$\d_m u$ is continuous for the product topology, with $\cP_2(\bS)$ equipped with the 2-Wasserstein distance, and has quadratic growth in $x \in \dbR^d$, locally uniformly in $m \in \cP_2(\bS)$, so as to guarantee integrability in the last expression. As in \cite{TTZ}, we denote
\bea
\label{dmu1}
\d_m u_i(t,m,x) := \d_m u(t,m,x,i) 
&\mbox{for}&
i \in \{0,1\},\qq
D_I u := \d_m u_1 - \d_m u_0,
\eea
and we introduce the measure flow generator of $X$
\bea
\label{cLx}
\left.\ba{c}
\dis \dbL u(t,m) := \pa_t u(t,m) + \int_{\dbR^d}  \cL_x\d_m u_1(t,m,x) m(dx,1),\\
\mbox{where}\q \cL_x\d_m u_1 := b\cd \pa_x \d_m u_1 + {1\over 2} \si^2 : \pa_{xx}^2  \d_m u_1.
\ea\right.
\eea
 Throughout this paper, we denote by 
 \beaa
 \mathbf{Q}_t:=[t,T)\times\cP_2(\mathbf{S}), 
 &\mbox{and}&
 \overline{\mathbf{Q}}_t:=[t,T]\times\cP_2(\mathbf{S}),
 ~~t\in[0,T).
 \eeaa
 
\begin{defn}\label{C12S}
Let $C^{1,2}_2(\overline{\mathbf{Q}}_t)$ be the set of functions $u: \overline{\mathbf{Q}}_t \to \dbR$ s.t. 
\\
$~\hspace{5mm}\bullet$ $\pa_t u, \d_m u, \pa_x \d_m u_1, \pa_{xx}^2\d_m u_1$ exist, and are continuous in all variables, 
\\
$~\hspace{5mm}\bullet$ $\pa_{xx}^2\d_m u_1$ is bounded in $x$, locally uniformly in $(t, m)$.
\end{defn}

The following Itô's formula is due to \cite[\S 3]{TTZ}:  for any $u\in C^{1,2}_2(\overline{\mathbf{Q}}_0)$ and $\dbP \in \cP(0,m)$:
 \bea
 \label{JumpFunIto}
 \left.\ba{lll}
\dis  u(T,m_{T-}) = u(0,m) + \int_0^T \dbL u(s,m_s)ds \\ 
\dis\qq  + \sum_{s \in J_{[0, T)}(\mathbf{m})} [u(s,m_s)-u(s, m_{s-})] + \dbE^\dbP\Big[\int_{J^c_{[0, T)}(\mathbf{m})} D_I u(s,m_s,X_s)dI_s \Big],
\ea\right.
\eea
where $\bm := \{m_s:= \dbP_{Y_s}\}_{s\in [-1, T]}$, $J_\dbT(\bm) := \{s \in \dbT: m_s \neq m_{s-}\}$, for all subset $\dbT\subset [0, T]$, and $J^c_\dbT(\bm)$ its complement set in $\mathbb{T}$.

\subsection{The dynamic programming equation}

Given two probability measures $m, m'\in\cP_2(\mathbf{S})$, we say that $m' \preceq m$ if $m'(\cdot,1)$ is absolutely continuous w.r.t.\ $m(\cdot, 1)$ with density bounded by $1$, i.e.
\bea\label{order}
m'(dx, 1) = p(x) m(dx, 1), \ \mbox{and} \ m'(dx,0) = [1-p(x)]m(dx,1) + m(dx, 0), 
\eea
for some measurable $p : \dbR^d \rightarrow [0,1]$. In other words, $m'(dx,1)$ is obtained from $m$ by randomly stopping a proportion $1-p(x)$ of the surviving particles. In our context, $m_{t-} = \dbP_{(X_t,I_{t-})}$ and $m_t = \dbP_{(X_t, I_t)}$, with $\dbP \in \cP(t,m)$, so that $m_t \preceq m_{t^-}$ with conditional transition probability $p(x) = p_t(x) := \dbP(I_t = 1 \mid X_t = x, I_{t-} = 1)$.

The following property (proved in Appendix \ref{appendix}) will be often used in this paper:

\begin{lem}\label{lem:compact-smallest}
For an arbitrary $m \in \cP_2(\bS)$:

\no{\rm (i)} the set $\{m' : m' \preceq m\}$ is compact,

\no{\rm (ii)} any compact subset $\cK(m)\subset\{m' : m' \preceq m \}$ has a smallest element for $\preceq$, i.e., there exists $\bar m \in \cK(m)$ such that for all $m' \in \cK(m)$, we have
$ m' \preceq \bar m$ implies that $m' = \bar m$. 
\end{lem}

The dynamic programming equation corresponding to our mean field optimal stopping problem is the infinitesimal counterpart of the dynamic programming principle \eqref{weakDPP}, and is defined by
\bea\label{obstacle}
\min\Big\{\min_{m' \in C_u(t,m)} \big[-(\dbL u + F)(t,m')\big],~ (\dbD_Iu)_* (t,m) \Big\} = 0, ~ (t,m) \in\mathbf{Q}_0,
\eea
with boundary condition $u|_{t=T} = g$. Here  the function $(\dbD_Iu)_*$ is the LSC envelope of 
 $$
 \dbD_Iu:(t,m) \longmapsto \inf_{x\in\mbox{\tiny Supp}(m(\cd,1))}D_I u(t,m,x), 
 $$
which is upper semicontinuous, but not continuous, in general,  and the set 
 \beaa
 C_u(t,m) 
 &:=& 
 \big\{m' \preceq m : u(t,m') \ge u(t,m)\big\},
 ~~(t,m) \in \mathbf{Q}_0,
 \eeaa
indicates the set of positions in $\mathbf{Q}_0$ which improve $u$ by stopping the corresponding particles.

For the purpose of the present paper, we note that this equation is slightly different from the obstacle equation introduced in our previous work \cite{TTZ}:

- if $u$ is a classical solution of \eqref{obstacle}, then it is nondecreasing for $\preceq$ (see \cite[Lemma 4.3]{TTZ}), and thus $C_u(t,m)$ is characterized by an equality, as in \cite{TTZ};

- despite the remaining differences, the two equations define the same solution, but this does not seem to have an immediate proof; we emphasize however that the equivalence between these two equations is a direct consequence of our uniqueness result in \cite[Theorem 4.5]{TTZ}, and the comparison Theorem \ref{thm-comparison} below.

Our objective in this paper is to develop a notion of viscosity solution for this equation which bypasses the strong regularity requirements of classical solutions. As usual, we start by introducing the notions of the sub- and supersolutions.
 
\begin{defn}
\label{defn-classical}
Let $u \in C^{1,2}_2(\overline{\mathbf{Q}}_0)$.
\\
{\rm (i)} $u$ is a classical supersolution of \reff{obstacle} if 
\bea
\label{classicalsuper}
\min\big\{ -(\dbL u + F) ,~\dbD_I u\big\}(t,m) \ge 0,\q
\forall (t, m)\in\mathbf{Q}_0. 
 \eea
 {\rm (ii)} $u$ is a classical subsolution of \reff{obstacle} if 
\bea
\label{classicalsub}
\min\big\{ -(\dbL u + F) , ~(\dbD_I u)_*\big\}(t,m) \le 0,
\q \forall (t,m)\in\mathbf{Q}_0~\mbox{s.t.}~C_u(t, m)  = \{m\}. 
 \eea
{\rm (iii)} $u$ is a classical solution of \eqref{obstacle} if it is a classical supersolution and subsolution.
\end{defn}

\section{Viscosity solutions }
\label{sect-viscosity}

\subsection{Definition and consistency}

For $\d > 0$ and $(t,m) \in \mathbf{Q}_0$, we introduce the neighborhood
$$ \cN_\d (t,m) := \big\{(s, \tilde m) :  s \in [t, t+\d], \dbP \in \cP(t,m), \ \tilde m \in \{ \dbP_{Y_{s-}}, \dbP_{Y_s} \} \big\}. $$
Note that, as the closure of a càdlàg $\cP_2(\bS)$-valued graph, $\cN_\d(t,m)$ is compact, by the compactnesses of $[t,t+\d]$, $\cP(t,m)$ and $\{ (\dbP_{Y_{s-}}, \dbP_{Y_s}) \}_{s \in [t,t+\d]}$ for any $\dbP \in \cP(t,m)$. 

\begin{defn}\label{NdUSC}
Let $u : \ol \bQ_0 \longrightarrow \dbR$. We say that $u$ is $\cN$-USC (respectively $\cN$-LSC) if
$$ u(t,m) \ge \underset{(s,\tilde m) \rightarrow (t,m)}{\lim \sup} u(s,\tilde m) \q \mbox{(respectively $u(t,m) \le \underset{(s,\tilde m) \rightarrow (t,m)}{\lim \inf} u(s,\tilde m)$)} \qq \mbox{for all $(t,m) \in \ol \bQ_0$}, $$
where the limits are sequences $(t_n, m_n) \rightarrow (t,m)$, with $(t_n,m_n) \in \cN_{T-t}(t,m)$.
\end{defn}

Note that the standard $\cW_2$-semicontinuity implies the $\cN$-semicontinuity. For a  locally bounded function $u: \overline{\mathbf{Q}}_0\longrightarrow\dbR$, we introduce its $\cN$-LSC and $\cN$-USC envelopes relatively to $\cP(t,m)$, $u_*$ and $u^*$ respectively:
\beaa
u_*(t,m) := \underset{(s, \tilde m) \to (t,m)}{\lim \inf} u(s,\tilde m), \q u^*(t, m) := \underset{(s,\tilde m) \to (t,m)}{\lim \sup} u(s,\tilde m), 
&\mbox{for all}&
(t,m) \in \overline{\mathbf{Q}}_0, 
\eeaa
where the limits are taken on all sequences $\{t_n,m_n\}_{n \ge 1}$ converging to $(t,m)$, with $(t_n,m_n) \in \cN_{T-t}(t, m)$ 
for all $n$. We then introduce the sets of test functions
\bea
\overline{\cA}u(t,m) &:=& \Big\{\f \in C^{1,2}_2(\overline{\mathbf{Q}}_t):  (\f-u_*)(t,m) = \max_{\cN_\d(t, m)}(\f-u_*) \ \mbox{for some $\d > 0$}\Big\} \nonumber, \\ 
\underline{\cA}u(t,m) &:=& \Big\{\f \in C^{1,2}_2(\overline{\mathbf{Q}}_t): (\f-u^*)(t,m) = \min_{\cN_\d(t, m)}(\f-u^*) \ \mbox{for some $\d > 0$}\Big\}. \nonumber
\eea

\begin{defn}
\label{defn-viscosity}
Let $u : \mathbf{Q}_0 \rightarrow \dbR$ be locally bounded.

\no{\rm (i)} $u$ is a viscosity supersolution of \reff{obstacle} if, for any  $(t, m)\in\mathbf{Q}_0$,
\bea
\label{super}
u_*(t,m) \ge u_*(t,m'),~ \forall  m' \preceq m,\q
\mbox{and}\q
-(\dbL \f + F) (t,m) \ge 0,
~\forall \f \in \overline{\cA}u(t,m). 
 \eea
 
 \no{\rm (ii)} $u$ is a viscosity subsolution of \reff{obstacle} if,  for any  $(t, m)\in\mathbf{Q}_0$ s.t.  $C_{u^*}(t, m)  = \{m\}$,
\bea
\label{sub}
\min\{ -(\dbL \f + F) , (\dbD_I \f)_*\}(t,m) \le 0,\q\forall~
 \f \in \underline{\cA}u(t,m).
 \eea

\no{\rm (iii)} $u$ is a viscosity solution of \eqref{obstacle} if it is a viscosity supersolution and subsolution.

\end{defn}

\begin{rem}\label{strict-visc} {\rm
Without loss of generality, we may assume that the maximum in the definition of $\ol \cA u(t,m)$ is strict. Indeed, for $\f \in \ol \cA u(t,m)$, we set
$$ \tilde \f(s, \tilde m) := \f(s,\tilde m) - (s-t)^2 - \big(\tilde m(\dbR^d,1) - m(\dbR^d, 1)\big)^2 \q \mbox{for all $(s, \tilde m) \in \bQ_t$.} $$
It is obvious that $\tilde \f\in C^{1,2}_2(\overline{\mathbf{Q}}_t)$. 
As $\tilde \f(s,\tilde m) = \tilde \f(t, m)$ if and only if $s = t$ and $\tilde m = m$ (since $\tilde m(\dbR^d,1) = m(\dbR^d, 1)$, and observing that in this case $\tilde m = \dbP_{Y_t} \preceq m$ for some $\dbP \in \cP(t,m)$), we deduce that $\tilde \f \in \ol \cA u(t,m)$ and the maximum is strict. Moreover, simple computations show that $\pa_t \tilde \f(t,m) = \pa_t \f(t,m)$ and $\dbL \tilde \f(t,m) = \dbL \f(t,m)$. An analogous statement holds for $\ul \cA u(t,m)$. 
\qed} \end{rem}

Our first result shows the consistency between classical and viscosity solutions.
\begin{thm}\label{thm-consistency}
Let $u \in C^{1,2}_2(\overline{\mathbf{Q}}_0)$. Then $u$ is a classical sub- (resp. super)solution of \eqref{obstacle} if, and only if, it is a viscosity sub- (resp. super)solution of \eqref{obstacle}.
\end{thm}
\proof
(i) Let $(t,m) \in \mathbf{Q}_0$. If $u$ is a viscosity super/subsolution, then given its smoothness we have $u \in \ol \cA u(t,m) \cap \ul \cA u(t,m)$, and we immediately deduce that $u$ is a classical super/subsolution. In particular, by \cite[Lemma 4.3]{TTZ}, $u$ is nondecreasing for $\preceq$ implies that $D_I u \ge 0$. 

\no(ii) Assume $u$ is a classical supersolution of \eqref{obstacle}. By \reff{classicalsuper} we see that $D_I u \ge 0$, then by \cite[Lemma 4.3]{TTZ} again we see that $u$ is nondecreasing for $\preceq$. Now let $\f \in \ol \cA u(t,m)$ with corresponding $\d$. 
 Introduce $\psi := \f-u$ and let 
$\bar\dbP \in \cP(t,m)$ be defined by \reff{barP} s.t. $X$ is unstopped under $\bar \dbP$. By definition of $\ol \cA u(t,m)$, we have
$ \psi(t,m) \ge \psi(s,\bar \dbP_{Y_s})$ for all $s \in [t,t+\d].$
Applying Itô's formula \eqref{JumpFunIto}, since the jump terms are equal to zero under $\bar \dbP$, we obtain:
$ - \frac{1}{\d}\int_t^{t+\d} \dbL \psi(s, \bar \dbP_{Y_s})ds \ge 0. $
Send $\d\to 0$, by the continuity of $s \mapsto \bar\dbP_{Y_s}$ we have $-\dbL \psi(t,m) \ge 0$, hence
$ -(\dbL \f + F) (t,m) \ge -(\dbL u + F) (t,m) \ge 0$
by the supersolution property of $u$.

Assume now that $u$ is a classical subsolution. Let $\f \in \ul \cA u(t,m)$ with corresponding $\d$, 
and assume that $(\dbD_I)_* \f(t,m) > 0$ and $C_u(t,m) = \{m\}$.  By definition of $\ul \cA u(t,m)$, we have
\bea
\label{consistency-sub}
 [\f-u](t,m)  \le [\f- u](s,\dbP_{Y_s}),\q \mbox{for all}~ s \in [t,t+\d], \dbP\in \cP(t,m). 
\eea
Set $s=t$ in \reff{consistency-sub}, then it follows from the arbitrariness of $\dbP\in \cP(t, m)$ that 
$ [\f-u](t,m)  \le [\f- u](t, m')$ for all $m' \preceq m$. 
Following the arguments of \cite[Lemma 4.3]{TTZ}, we deduce from above that $D_I [\f-u](t, m, \cd) \le 0$, 
and therefore $(\dbD_I u)_* (t,m) \ge (\dbD_I \f)_* (t,m) > 0$. The subsolution property of $u$ then implies that $-(\dbL u + F) (t,m) \le 0$. Using Itô's formula \eqref{JumpFunIto} under $\bar \dbP$ again on $[t, t+\d]$, we get from \reff{consistency-sub} that
$ -(\dbL \f + F) (t,m) \le -(\dbL u + F) (t,m) \le 0.$
\qed

\subsection{Some regularity results}
In this subsection, we present some regularity results which will be used in the rest of this section. Since our main focus is the viscosity properties, we postpone their proofs  to Appendix \ref{appendix}. 

\begin{lem}
\label{lem-localbound}
Under Assumption \ref{assum-bsig}, the value function $V$ is USC under $\cW_2$.
\end{lem}

\begin{thm}
\label{thm-reg}
{\rm (i)} Assume $f$ and $g$ are uniformly continuous in $(t, x, m)$, under $\cW_2$ for $m$, then $V$ is  continuous on $\overline{\bQ}_0$, under $\cW_2$ for $m$.

\no{\rm (ii)}  Assume further that $b, \si$ are uniformly Lipschitz continuous in $m$ under $\cW_1$, and $f, g$ are uniformly continuous in $m$ under $\cW_1$,   then $V$ is also continuous in $m$ under  $\cW_1$.  
\end{thm}

Even for the standard optimal stopping problems, one can hardly expect the value function to be smooth. We next establish a regularity result for the value function when $X$ is unstopped.  For $(t,m) \in \bQ_0$, let $\bar \dbP^{t, m} \in \cP(t,m)$ be as by \reff{barP}, and define
\bea
\label{XU}
 U(t, m) := g(\bar \dbP_{Y_T}^{t,m})  + \int_t^T F(r, \bar\dbP_{Y_r}^{t,m})dr.
\eea

\begin{lem}
\label{lem-Ureg}
For $\f=b, \si, f, g$,  assume $\f$ is continuous in $t$ and $\pa_x\f, \d_m \f, \pa_x \d_m \f, \pa_{xx}^2\d_m\f$ exist and are continuous and bounded  and that, for $\f = b, \si$, all ther derivatives of $\f$ are Lipschitz up to order 2.  Then $U \in C^{1,2}(\overline\bQ_0)$ with bounded  $\pa_x \d_m U, \pa_{xx}\d_m U$ and in particular $U \in C_2^{1,2}(\overline\bQ_0)$. Moreover, if $b, \si, f, g$ are  uniformly Lipschitz continuous in $m$ under $\cW_1$ with a Lipschitz constant $L$, then $U$ is uniformly Lipschitz continuous in $m$ under $\cW_1$ with a Lipschitz constant $C_L$.
\end{lem}

Finally, we introduce a smooth mollifier for functions on the Wasserstein space. 

  \begin{lem}\label{lem-mol}

{\rm (i)} Let $U : \cP_2(\bS) \rightarrow \dbR$ be continuous. There exists $\{U_n\}_{n \ge 1}$ in $C^\infty(\cP_2(\bS))$ such that $\dis\lim_{n\to\infty} \sup_{m\in \cM}|U_n(m) - U(m)| =0$ for any compact set  $\cM\subset \cP_2(\bS)$;\\
  {\rm (ii)} Let $U : \cP_1(\bS) \rightarrow \dbR$ be continuous under $\cW_1$. There exists $\{U_n\}_{n \ge 1}$ in $C^\infty(\cP_2(\bS))\cap C^0(\cP_1(\bS))$ such that  $\dis\lim_{n\to\infty} \sup_{m\in \cM}|U_n(m) - U(m)| =0$ for any compact  set $\cM\subset \cP_1(\bS)$;\\
{\rm (iii)} Assume further that $U$ is Lipschitz continuous under $\cW_1$, then we may choose $\{U_n\}_{n \ge 1}$ to be Lipschitz continuous under $\cW_1$, uniformly in $n$. 
\end{lem}
The mollifier is adopted from Mou \& Zhang \cite{MZ}. Note that the extension of the state space from $\dbR^d$ in \cite{MZ} to $\bS$ here is straightforward.
We remark that if $U$ is Lipschitz continuous under $\cW_2$,  in general the Lipschitz continuity of $U_n$ under $\cW_2$ is not uniform in $n$.

 \subsection{Viscosity property}
We first need a simple lemma whose proof is postponed to Appendix  \ref{appendix}. 
 
 \begin{lem}\label{lem-lusin}
{\rm (i)} Let $v : \cP_2(\bS) \longrightarrow \dbR$ be $\cN$-LSC, and $m \in \cP_2(\bS)$ s.t. $v(m) \ge v(m')$ for all $m' \preceq m$ with continuous conditional transition probability. Then $v(m) \ge v(m')$ for all $m' \preceq m$.

\no{\rm (ii)} Let $\f \in C^0(\ol \bQ_0, \dbR)$ admit a continuous linear functional derivative $\d_m \f$. Assume we have $(\dbD_I \f)_*(t,m) > 0$ for some $(t,m) \in \bQ_0$. Then $\f$ is nondecreasing for $\preceq$ in a neighborhood of $(t,m)$.
\end{lem}

\begin{thm}
\label{thm-existence}
The value function $V$  is a viscosity solution of \eqref{obstacle}.
\end{thm}
\proof First, by Lemma \ref{lem-localbound},  $V$ inherits the local boundedness of $g$.

 \no(i) We first verify the viscosity supersolution property. Fix $(t, m)$ and $\f\in \overline{\cA}V(t,m)$. We may assume w.l.o.g. that $[V_*-\f](t,m) = 0$. Let $\d > 0$ and $(t_n, m_n)_{n \ge 1} \in \cN_\d(t,m)$ converging to $(t,m)$ s.t. $V(t_n, m_n) \underset{n \to \infty}{\longrightarrow} V_*(t,m)$, and denote $\eta_n := [V-\f](t_n,m_n) \ge 0$, as $V\ge V_*$. Thus, we have $\eta_n \underset{n \to \infty}{\longrightarrow} 0$. By the DPP \eqref{weakDPP}, we have 
\beaa
\eta_n + \f(t_n, m_n)
&=&
V(t_n, m_n) 
\ge \int_{t_n}^{s_n} F(r, \bar \dbP_{Y_r}^{m_n})dr + V(s_n, \bar \dbP_{Y_{s_n}}^{m_n}) \\
&\ge& \int_{t_n}^{s_n} F(r, \bar \dbP_{Y_r}^{m_n})dr + V_*(s_n, \bar \dbP_{Y_{s_n}}^{m_n}) \ge \int_{t_n}^{s_n} F(r, \bar \dbP_{Y_r}^{m_n})dr + \f(s_n, \bar \dbP_{Y_{s_n}}^{m_n}), 
\eeaa
where $\bar \dbP^{m_n}:=\bar \dbP^{t_n, m_n} \in \cP(t_n, m_n)$ is defined by \reff{barP} such that $X$ is unstopped, and $s_n := t_n + h_n$ with $h_n := \sqrt{\eta_n} \vee n^{-1}$.
Thus, by Itô's formula, the above gives
$ h_n + \frac{1}{h_n}\int_{t_n}^{s_n} -(\dbL \f + F)(r, \bar \dbP_{Y_{r}}^{m_n})dr \ge 0.$
Send $n\to \infty$, since $h_n{\longrightarrow} 0$, we obtain  $-(\dbL \f + F)(t,m) \ge 0$. 

We now prove the remaining part of the supersolution property. Let $m' \preceq m$ with transition probability $p$. By Lemma \ref{lem-lusin} (i), we may assume without loss of generality that $p$ is continuous. For all $n \ge 1$, define $m_n' \preceq m_n$ as the measure obtained from $m_n$ by applying the same $p$. Given the continuity of $p$ and the compactness $\cN_\d(t,m)$, we see by \eqref{order} that $\cW_2(m_n', m') \underset{n \to \infty}{\longrightarrow} 0$. Let $\bar \dbP^{m_n, m_n'} \in \cP(t_n, m_n)$ be s.t. $\bar \dbP_{Y_{t_n}}^{m_n, m_n'} = m_n'$, and $I_s = I_{t_n}$, $\bar \dbP^{m_n, m_n'}$-a.s. for all $s \ge t_n$.
By \eqref{weakDPP}, 
 \bea\label{superexistence}
 V(t_n, m_n) \ge \int_{t_n}^{s} F(r, \bar \dbP_{Y_r}^{m_n, m_n'})dr + V\big(s, \bar \dbP_{Y_{s}}^{m_n, m_n'}\big), 
 &\mbox{for all}&
 s \ge t_n
 ~~\mbox{and}~~
 n \ge 1. 
 \eea
Take $s=t_n$ and $\lim \inf_{n\to\infty}$ in \reff{superexistence}, we obtain $ V_*(t,m) \ge V_*(t,m')$ as $V(t_n, m_n) \to V_*(t,m)$.

\no(ii) We next verify the viscosity subsolution property.   Let $(t,m)$ and $\f \in \ul \cA V(t,m)$ be s.t. $C_{V^*}(t,m) = \{m\}$ and $(\dbD_I\f)_* (t,m) > 0$. We may assume w.l.o.g. that $[V^*-\f](t,m) = 0$.  Let $\d > 0$ and $(t_n, m_n)_{n \ge 1} \in \cN_\d(t,m)$ converging to $(t,m)$ such that $V(t_n, m_n) \underset{n \to \infty}{\longrightarrow} V^*(t,m)$, and denote $-\eta_n := [V-\f](t_n,m_n) \le 0$.  Thus $\eta_n \underset{n \to \infty}{\longrightarrow} 0$. For $n\ge 1$, since $g$ is USC and $\cP(t_n, m_n)$ is compact, there exists  $\dbP^{n,*} \in \cP(t_n, m_n)$ s.t. $V(t_n, m_n) = \int_{t_n}^{T} F(r, \dbP_{Y_r}^{n,*})dr + g(\dbP_{Y_T}^{n,*})$. By DPP, we have
$$ V(t_n, m_n) \ge \int_{t_n}^{s_n} F(r, \dbP_{Y_r}^{n,*})dr + V(s_n, \dbP_{Y_{s_n}}^{n,*}) \ge \int_{t_n}^{T} F(r, \dbP_{Y_r}^{n,*})dr + g(\dbP_{Y_T}^{n,*}),$$
where $s_n := t_n + h_n$ with $h_n := \sqrt{\eta_n}\vee n^{-1}$, and thus,
\bea\label{optDPP}
V(t_n, m_n) = \int_{t_n}^{s_n} F(r, \dbP_{Y_r}^{n,*})dr + V(s_n, \dbP_{Y_{s_n}}^{n,*}).
\eea 
Note that $\dbP^{n,*} \in \cP(t_n, m_n) \subset \cP(t,m)$ for all $n$, and $\cP(t,m)$  is compact, we may extract a subsequence (still denoted the same) s.t. $\dbP^{n,*}  \underset{n \to \infty}{\longrightarrow}  \dbP^*$, for some $\dbP^* \in \cP(t,m)$. As the trajectories $r \mapsto \dbP_{Y_{r}}^{n,*}$ are càdlàg and $s_n \downarrow t$, this implies $\cW_2(\dbP_{Y_{s_n}}^{n,*}, m^*) \underset{n \to \infty}{\longrightarrow} 0$, where $m^* := \dbP_{Y_t}^* \preceq m$ as $\dbP^* \in \cP(t,m)$. 
Thus, taking the $\lim \sup_{n \to \infty}$ in \eqref{optDPP} and recalling $V(t_n, m_n) \to V^*(t,m)$, we have 
$ V^*(t,m) \le V^*(t,m^*).$
 As $C_{V^*}(t,m) = \{m\}$, we obtain $m^* = m$. Moreover, \eqref{optDPP} also implies that
\bea\label{subexistence}
 -\eta_n + \f(t_n, m_n) \le \int_{t_n}^{s_n} F(r, \dbP_{Y_r}^{n,*})dr + \f(s_n, \dbP_{Y_{s_n-}}^{n,*}) \q \mbox{for all $n \ge 1$.}
 \eea
Let $\cB_{\cW_2}(m,\d)$ denote the $\cW_2$ ball centered in $m$, with radius $\d$. By Lemma \ref{lem-lusin} (ii), the fact that $(\dbD_I \f)_*(t,m) > 0$ implies that $\f$ is (strictly) increasing for $\preceq$ on $[t,t+\d) \times \cB_{\cW_2}(m,\d)$, for a possibly smaller $\d>0$. By convergence to $(t,m)$, we have $\{\dbP_{Y_{r}}^{n,*}, t\le r\le s_n\}  \subset \cB_{\cW_2}(m,\d)$ for $n$ large. Then $D_I \f(r,\dbP_{Y_{r-}}^{n,*},\cd) \ge 0$ and $\f(r,\dbP_{Y_{r-}}^{n,*})\ge \f(r,\dbP_{Y_{r}}^{n,*})$ for $t\le r\le s_n$. Using the fact that the trajectories are càdlàg, by applying  It\^{o}'s formula on \eqref{subexistence} we obtain $-(\dbL \f + F) (t,m) \le 0$. \qed

\subsection{Stability}

\begin{thm}
Let $\{F_\e\}_{\e > 0}$ be a family of functions from $\ol \bQ_0$ to $\dbR$ such that $F_\e \underset{\e \to 0}{\longrightarrow} F$ uniformly on compact subsets of $\ol \bQ_0$, and let $\{u_\e\}_{\e > 0}$ and $\{v_\e\}_{\e > 0}$ be two families of viscosity subsolutions and supersolutions of \eqref{obstacle} with $F_\e$ instead of $F$, respectively. Assume that the following relaxed semilimits are finite
\beaa
\ol u(t,m) := \underset{(\e,s,\tilde m) \rightarrow (0,t,m)}{\lim \sup} u_\e(s,\tilde m),
~\mbox{and}~
\ul v(t,m) := \underset{(\e,s,\tilde m) \rightarrow (0,t,m)}{\lim \inf} v_\e(s,\tilde m),
&(t,m) \in \overline{\mathbf{Q}}_0,&
\eeaa
where the limits are sequences $(\e_n, t_n, m_n) \rightarrow (0,t,m)$, with $(t_n,m_n) \in \cN_{T-t}(t,m)$. Then $\ol u$ (resp. $\ul v$) is a $\cN$--USC (resp. $\cN$--LSC) viscosity subsolution (resp. supersolution) of \eqref{obstacle}.  
\end{thm}
\proof
(i) We prove the stability of the supersolution first. Observe that we may assume without loss of generality that $v_\eps$ is $\cN$--LSC as $
\ul v(t,m) = \underset{(\e,s,\tilde m) \rightarrow (0,t,m)}{\lim \inf} (v_\e)_*(s,\tilde m)$. Also note that $\ul v$ is clearly $\cN$-LSC in the sense of Definition \ref{NdUSC}.

Fix $(t,m) \in {\mathbf{Q}}_0$, and $\f \in \ol \cA \ul v(t,m)$ with corresponding $\d$, and s.t. $(t,m)$ is a strict maximizer of $\f - \ul v$ on $\cN_\d(t,m)$, see Remark \ref{strict-visc}. By definition, there exists a sequence $(\e_n, t_n, m_n)\to (0, t, m)$ s.t. $v_{\e_n}(t_n,m_n) \rightarrow \ul v(t,m)$. Note that $(t_n, m_n) \in \cN_\d(t,m)$ for all $n$ large, then we can find $\d' < \d$ s.t. $\cN_{\d'}(t_n,m_n) \subset \cN_\d(t,m)$. Let $(\hat t_n, \hat m_n)$ be a maximizer of $\f -v_{\e_n}$ on $\cN_{\d'}(t_n,m_n)$.  We first note that
\bea\label{conv-stable}
(\hat t_n, \hat m_n) \underset{n \rightarrow \infty}{\longrightarrow} (t,m).
\eea
Indeed,  $(\hat t_n, \hat m_n) \in \cN_{\d'}(t_n,m_n) \subset \cN_\d(t,m)$ for all $n$. Thus, by compactness, there exists a subsequence (still named $\hat m_n$) converging to some $(\hat t, \hat m) \in \cN_\d(t,m)$. Observing that 
\beaa
[\f - \ul v](t,m) &=& \underset{n \rightarrow \infty}{\lim}[\f - v_{\e_n}](t_n,m_n) \le \underset{n \rightarrow \infty}{\lim \inf}[\f - v_{\e_n}](\hat t_n, \hat m_n)\\
&\le&  \underset{n \rightarrow \infty}{\lim \sup}[\f - v_{\e_n}](\hat t_n, \hat m_n) \le [\f - \ul v](\hat t, \hat m),
\eeaa
we conclude from the fact that $(t,m)$ is a strict maximizer of $\f - \ul v$ on $\cN_\d(t,m)$ that $(\hat t, \hat m) = (t,m)$, and thus \eqref{conv-stable} holds true. Then, given that $(t_n, m_n)$ and $(\hat t_n, \hat m_n)$ have the same limit, we have $\cN_{\d''}(\hat t_n, \hat m_n) \subset \cN_{\d'}(t_n, m_n)$ for some $\d'' < \d'$ and $n$ large enough. Then, as $(\hat t_n, \hat m_n)$ is also a maximizer on $\cN_{\d''}(\hat t_n, \hat m_n)$, the supersolution property implies
$ -(\dbL \f + F_{\e_n}) (\hat t_n, \hat m_n) \ge 0,$ for $n$ large enough, 
and we derive the first part of the supersolution property of $\ul v$ by sending $n \rightarrow \infty$. 

We now prove that $\ul v$ is increasing for $\preceq$. By Lemma \ref{lem-lusin} (i), it suffices to prove that $\ul v(t,m) \ge \ul v(t,m')$ for a given $m' \preceq m$ with continuous conditional transition probability $p$. We define for all $n$ the measure $m_n' \preceq m_n$, obtained from $m_n$ by applying $p$. As $\cW_2(m_n, m) \underset{n \to \infty}{\longrightarrow} 0$ and $p$ is continuous, similarly to the proof of Theorem \ref{thm-existence}, we see that $\cW_2(m_n', m') \underset{n \to \infty}{\longrightarrow} 0$. Moreover, by the supersolution property of $v_{\e_n}$, we have 
$ v_{\e_n}(t_n, m_n) \ge v_{\e_n}(t_n, m_n')$, for all $n \ge 1$,
 and we conclude by taking the $\lim \inf$ that $\ul v(t,m) \ge \ul v(t,m')$, as the l.h.s. of the inequality converges. 
 
\no(ii) We now prove the stability of the subsolution. Similarly to (i), we may assume that $\{u_\e\}_{\e > 0}$ is a family of $\cN$--USC viscosity subsolutions of \eqref{obstacle}, and observe that $\ol u$ is clearly $\cN$-USC.
 Let $(t,m)$ and $\f \in \ul \cA \ol u(t,m)$ be such that $(t,m)$ is a strict local minimizer of $\varphi-\overline{u}$. Assume that $C_{\ol u}(t,m) = \{m\}$ and $(\dbD_I \f)_*(t,m) > 0$.  Following the same argument as in the previous step, replacing maximizers with minimizers, we may construct $(\hat t_n, \hat m_n)$,  converging to some $(\hat t,\hat m)$, and satisfying the inequalities
$$[\f - \ol u](t,m)  \ge \underset{n \rightarrow \infty}{\lim \sup} \ [\f - u_{\e_n}](\hat t_n, \hat m_n) \ge \underset{n \rightarrow \infty}{\lim \inf} \ [\f - u_{\e_n}](\hat t_n, \hat m_n) \ge [\f - \ol u](\hat t, \hat m). $$
By the strict minimum property of $(t,m)$, this again implies that $(\hat t,\hat m)=(t,m)$, and $\underset{n \rightarrow \infty}{\lim} u_{\e_n}(\hat t_n, \hat m_n) = \ol u(t,m)$.
 By Lemma \ref{lem:compact-smallest}, we may  now take $m_n^* \in \underset{\preceq}{\mathrm{argmin}} \ C_{u_{\e_n}}(\hat t_n, \hat m_n)$. By compactness, there is a subsequence $\{m^*_n\}_{n \ge 1}$ converging to some $m^*$. As $u_{\e_n}(\hat t_n, \hat m_n) \le u_{\e_n}(\hat t_n, m_n^*)$ for all $n$, by definition of $C_{u_{\e_n}}(\hat t_n, \hat m_n)$, taking the lim sup implies $\ol u(t,m) \le \ol u(t,m^*)$, hence $m^* = m$ as $C_{\ol u}(t,m) = \{m\}$. As $(\dbD_I \f)_*(t,m) > 0$, $\hat m_n$ and $m_n^*$ are both in a neighborhood where $\f$ is strictly increasing for $n$ large enough, and thus
 $ [\f-u_{\e_n}](\hat t_n, \hat m_n) \ge [\f-u_{\e_n}](\hat t_n, m_n^*)$,
 which implies equality by definition of $(\hat t_n, \hat m_n)$ and the fact that $(\hat t_n, m_n^*) \in \cN_{\d''}(\hat t_n, \hat m_n)$.
 Then $\f \in \ul \cA u_{\e_n}(\hat t_n, m_n^*)$. As $C_{u_{\e_n}}(\hat t_n, m_n^*) = \{m_n^*\}$, the viscosity subsolution implies $-(\dbL \f + F_{\e_n}) (\hat t_n, m_n^*) \le 0$ for $n$ large enough, and we conclude by letting $n \longrightarrow \infty$. \qed

\begin{rem}{\rm
A natural extension of the stability result is to allow the perturbation of $b$ and $\si$. However, this would change the definition of $\cP(t,m)$ in \eqref{martingalepb}, and therefore our viscosity neighborhoods $\cN_\d(t,m)$. Although we expect the stability property to remain true, this would require to extend $\cP(t,m)$ in some sense, which would go beyond the scope of the present paper. 
}
\end{rem}

\subsection{Comparison}

\begin{thm}
\label{thm-comparison}
{\rm (i)} Let $u$ be an $\cN$-USC viscosity subsolution of \reff{obstacle} satisfying $u|_{t=T} \le g$. Assume further that $f$ is uniformly continuous in $(t,x,m)$ under $\cW_2$. Then $u\le V$.
\\
{\rm (ii)} Let $v$ be a $\cN$-LSC viscosity supersolution of \reff{obstacle} satisfying $v|_{t=T} \ge g$.  
Assume further that $b$, $\si$, $f$, $g$ can be extended to $\cP_1(\bS)$ under $\cW_1$ continuously; $b$ is uniformly Lipschitz continuous in $(x, m)$ under $\cW_1$; and $\si$ has the regularity required in Lemma \ref{lem-Ureg}. 
Then $v \ge V$. 
\end{thm}
\proof 
(i) We first compare $V$ and $u$. Assume by contradiction that  
$
u(t, m) > V(t,m) 
$ for some $(t, m)$. Then, for $\e>0$ small enough, 
\bea
\label{utTe}
u(t, m) - \f_\e(t,m) > \sup_{\dbP\in \cP(t, m)} \Big\{ \int_t^T F(r, \dbP_{Y_r})dr +  g(\dbP_{Y_T}) \Big\}, 
\eea
where $\f_\e(s, \tilde m) := \e \big[(T-t)+ m(\dbR^d,1)\big]$. Let $(t^*, \dbP^*)$ be s.t.
\bea
\label{c*}
\left.\ba{c}
\dis (u-\f_\e)(t^*,m^*) + \int_t^{t^*} \!\!\! F(r, \dbP_{Y_r}^*)dr 
\dis = \!\!\! \max_{\tiny \begin{array}{c} (s, m, \dbP)\in \cN_{T-t}(t,m) \times \cP(t,m) : \\ m \in \{ \dbP_{Y_{s-}}, \dbP_{Y_s} \} \end{array}} \!\!\! \Big\{(u-\f_\e)(s, m) +  \int_t^s \!\!\! F(r, \dbP_{Y_r})dr \Big\}.
\ea\right.
\eea
where $m^*$ is the optimal argument in $\{ \dbP_{Y_{t^*-}}^*, \dbP_{Y_{t^*}}^* \}$. 
Clearly $t^*<T$. Indeed, if $t^* = T$, then $(T,m^*) \in \cN_{T-t}(t,m)$, and by \reff{c*} and \reff{utTe} we have
\beaa
&&u(T, m^*) - \e m^*(\dbR^d,1) + \int_t^T F(r, \dbP_{Y_r}^*)dr  \ge (u-\f_\e)(t,m) \\
&&> \sup_{\dbP\in \cP(t, m)} \Big\{ \int_t^T F(r, \dbP_{Y_r})dr  + g(\dbP_{Y_T}) \Big\} \ge  \int_t^T F(r, \dbP_{Y_r}^*)dr + u(T, m^*), 
\eeaa
as $u(T, \cdot) \le g$. This is the desired contradiction. 
Moreover, by Lemma \ref{lem:compact-smallest}, we may choose $m^*$ to be the smallest one which keeps the same value $(u-\f_\e)(t^*, m^*)$. Note that this change is only at $t^*$, and thus has no impact on the value of  $\int_t^{t^*} F(r, \dbP_{Y_r}^*)dr$. Then
\bea
\label{tm*}
\left.\ba{c}
(u-\f_\e)(t^*, m^*)  > (u-\f_\e)(t^*, m'),\q \mbox{for all}~ m^*\neq m'\preceq m^*.
\ea\right.
\eea
Furthermore, we note that, since $m \mapsto m(\dbR^d,1)$ is increasing, by  \reff{tm*} actually we have
\bea
\label{Cu}
u(t^*, m^*) > u(t^*, m'),\q \mbox{for all} \ m^*\neq m'\preceq m^*,\q\mbox{namely}\q C_{u}(t^*, m^*) = \{m^*\}.
\eea 

Next, let $f^+, f^-$ denote the positive and negative part of $f$, respectively, and $\rho_0$ the modulus of continuity function of $f$. Introduce: \beaa
\underline f^+(s,x, \tilde m) :=  f^+(s,x, \tilde m) - \rho_0\big(\big|\tilde m(\dbR^d, 1) - m^*(\dbR^d,1)\big|^{1\over 2}\big);\\
\overline f^-(s,x, \tilde m) :=  f^-(s,x, \tilde m) +  \rho_0\big(\big|\tilde m(\dbR^d, 1) - m^*(\dbR^d,1)\big|^{1\over 2}\big).
\eeaa
It is clear that $\underline f^+, \overline f^-$ are also uniformly continuous in $(s,x,\tilde m)$ (under $\cW_2$).  For $\e>0$, by Lemma \ref{lem-mol} (i) let $\underline f^+_n, \overline f^-_n$ be a smooth mollifier (under $\cW_2$) such that 
\beaa
|\underline f^+_n-\underline f^+|\le {\e\over 6},\q |\overline f^-_n-\overline f^-|\le {\e\over 6},\q\mbox{on}\q \cP(t^*, m^*).
\eeaa
Then,  for all $(s, \tilde m) \in \cN_{T-t^*}(t^*, m^*)$ with corresponding $\dbP \in \cP(t^*, m^*)$, and $t^*\le r< s$, considering the case $\tilde m = \dbP_{Y_{s-}}$, we have
\beaa
&&\dis \dbE^\dbP\Big[\Big|f^+(s, X_{s}, \dbP_{(X_{s}, I_r)}) I_{s-}- f^+(s, X_{s}, \tilde m) I_{s-}\Big|\Big] \le  \rho_0\Big(\cW_2\big( \dbP_{(X_{s}, I_r)}, \dbP_{(X_s, I_{s-})}\big)\Big)\\
&&\dis\le  \rho_0\Big(\sqrt{\dbE^\dbP[|I_r - I_{s-}|^2]}\Big)\le   \rho_0\Big(\sqrt{\dbE^\dbP[|I_{t^*} - I_{s-}|^2]}\Big) = \rho_0\big(\big|\tilde m(\dbR^d, 1) - m^*(\dbR^d,1)\big|^{1\over 2}\big).
\eeaa
Then
 $\dbE^\dbP\Big[f^+(s, X_{s}, \dbP_{(X_{s}, I_r})) I_{s-})\Big] \ge   \dbE^\dbP\Big[\underline f^+(s, X_{s}, \tilde m) I_{s-}\Big]$
and, similarly,
$
 \dbE^\dbP\Big[f^-(t^*, X_{t^*}, \dbP_{(X_s, I_r)}) I_{t^*}\Big] \le   \dbE^\dbP\Big[\overline f^-(t^*, X_{t^*}, \tilde m) I_{t^*} \Big].
 $
Thus, by \reff{F} and the regularity of $f$, we have
\beaa
&&\dis F(r, \dbP_{Y_r})= \dbE^\dbP\Big[f(r, X_r, \dbP_{Y_r}) I_r\Big]=\dbE^\dbP\Big[f^+(r, X_r, \dbP_{Y_r}) I_r- f^-(r, X_r, \dbP_{Y_r}) I_r\Big]\\
&&\dis \ge  \dbE^\dbP\Big[f^+(r, X_r, \dbP_{Y_r}) I_{s-}- f^-(r, X_r, \dbP_{Y_r}) I_{t^*}\Big]\\
&&\dis \ge  \dbE^\dbP\Big[f^+(s, X_s, \dbP_{(X_s, I_r})) I_{s-} - f^-(t^*, X_{t^*}, \dbP_{(X_s, I_r)}) I_{t^*}\Big]- \rho(s-t^*)\\
&&\dis \ge  \dbE^\dbP\Big[\underline f^+(s, X_s, \tilde m) I_{s-} - \overline f^-(t^*, X_{t^*}, \tilde m) I_{t^*}\Big]- \rho(s-t^*)\\
&&\dis \ge  \dbE^\dbP\Big[\underline f^+_n(s, X_s, \tilde m) I_{s-} - \overline f^-_n(t^*, X_{t^*}, \tilde m) I_{t^*}\Big] -{\e\over 3}- \rho(s-t^*),
\eeaa
for some modulus of continuity $\rho$ which can be chosen to be smooth on $(0,\infty)$. That is,
\beaa
 F(r, \dbP_{Y_r})\ge \int \underline f_n^+(s, x, \tilde m) i~ \tilde m(dx,di)  - \int \overline f_n^-(t^*, x, \tilde m) i \ m^*(dx,di)  -{\e\over 3}- \rho(s-t^*).
\eeaa
In the case $\tilde m = \dbP_{Y_s}$, following similar arguments we see the above still holds true. 
 Denote
\beaa
&&\dis \phi^n_\e(s, \tilde m) := \int f^+_n(s,x, \tilde m)i \ \tilde m(dx,di) - \int f^-_n(t^*, x,\tilde m)i \ m^*(dx,di),\\
&&\dis \psi^n_\e(s, \tilde m) :=  \f_\e(s, \tilde m) - (s- t^*)\big[ \phi^n_\e(s, \tilde m)  - {\e\over 3}  - \rho(s-t^*)\big],
\eeaa
which are obviously in $C_2^{1,2}(\overline\bQ_0)$.  Then, by \reff{c*},
\beaa
\dis &&\dis (u-\psi^n_\e)(t^*, m^*) =(u-\f_\e)(t^*, m^*) \ge (u-\f_\e)(s, \tilde m) + \int_{t^*}^s F(r, \dbP_{Y_r})dr\\
\dis &&\ge  (u-\f_\e)(s, \tilde m) + (s- t^*) \big[ \phi^n_\e(s, \tilde m)  - {\e\over 3}  - \rho(s-t^*)\big] = (u-\psi^n_\e)(s, \tilde m).
\eeaa
 Thus $\psi^n_\e \in \ul \cA u(t^*, m^*)$. Note that
 \beaa
&&\dis  [\dbL \psi^n_\e + F] (t^*, m^*) = \Big[\dbL \f_\e - \phi^n_\e + {\e\over 3} + F\Big](t^*, m^*) = -\e - \phi^n_\e(t^*, m^*) + {\e\over 3} + F(t^*, m^*)\\
 &&\dis\qq = -{2\e\over 3} + \int \big[\underline f^+_n(t^*, x, m^*) - \overline f^-_n(t^*, x, m^*)\big]i \ m^*(dx,di)  + F(t^*, m^*)\\
 &&\dis\qq \le -{\e\over 3} + \int \big[\underline f^+(t^*, x, m^*) - \overline f^-(t^*, x, m^*)\big]i \ m^*(dx,di)  + F(t^*, m^*)\\
 &&\dis\qq \le -{\e\over 3} + \int \big[f^+(t^*, x, m^*) - f^-(t^*, x, m^*)\big]i \ m^*(dx,di)   + F(t^*, m^*) = -{\e\over 3} <0;\\
 &&\dis (D_I \psi^n_\e)(s, \tilde m) = \e - (s- t^*) (D_I \phi^n_\e)(s, \tilde m),\q\mbox{and thus}\q (D_I \psi^n_\e)_*(t^*, m^*) = \e >0. 
 \eeaa
 Recall \reff{Cu}, this contradicts the viscosity subsolution property of $u$.

\medskip

\noindent (ii) We next compare $V$ and $v$. Fix $\e>0$. For each $n\ge 1$, denote $t_i := t^{(n)}_i:=\frac {iT} n$, $0 \le i \le n$. First, note that, for $(t,m) \in \mathbf{Q}_0$, it follows from the continuity of the coefficients that
\bea
\label{Vn}
\left.\ba{c}
\dis V(t, m) := \lim_{n\to\infty} V_n(t, m),\q \mbox{where}\q V_n(t, m) := \sup_{\dbP\in \cP_n(t, m)} \Big\{ \int_t^T F(r, \dbP_{Y_r})dr + g(\dbP_{Y_T}) \Big\},\\
\dis \cP_n(t, m):= \Big\{\dbP\in \cP(t, m): \mbox{$\t$ takes values in $\{t_1,\ldots,t_n\} \cap [t, T]$, $\dbP$-a.s.}\Big\}.
\ea\right. \nonumber
\eea

{\it Step 1:} We show that $(V_n - v)(t_{n-1}, \cdot) \le \frac \e n$. Assume to the contrary that there exists $m_{n-1}$ such that $(V_n-v)(t_{n-1}, m_{n-1}) > \frac \e n$. By the definition of $\cP_n(t_{n-1}, m_{n-1})$, we have $V_{n}(t, m)= \int_t^T F(r, \bar \dbP^{t,m}_{Y_r})dr + g(\bar \dbP_{Y_T}^{t,m})$, $t\in (t_{n-1}, T]$,  where $\bar \dbP^{t, m} \in \cP(t,m)$ is defined by \reff{barP}. 

Let $\d_1, \d_2>0$ be small numbers which will be specified later. Applying Lemma \ref{lem-mol} (ii), (iii), let $(g_k, f_k, b_k)$ be the smooth mollifier of $(g, f, b)$ (under $\cW_1$), where $b_k$ is also mollified in $(t,x)$ in a standard way, such that $\lVert g_k-g \rVert_\infty + \lVert f_k-f \rVert_\infty  \le \d_1$, $\lVert b_k - b \rVert_\infty \le \d_2$, and $g_k$ is Lipschitz continuous under $\cW_1$ with a Lipschitz constant $L_k$ depending on $k$, and $b_k$ is uniformly Lipschitz continuous in $(x, m)$ under $\cW_1$ with a Lipschitz constant $L$ independent of $k$. By otherwise choosing a larger $L$ we assume $\si$ is also uniformly Lipschitz continuous in $(x, m)$ under $\cW_1$ with Lipschitz constant $L$. Let $U^{k_1, k_2}$ be defined by \reff{XU} corresponding to $(b_{k_2}, \si, g_{k_1}, f_{k_1})$. Then, by Lemma \ref{lem-Ureg}, 
\bea\label{eqdiff}
\pa_tU^{k_1,k_2}(t,m) + \int_{\dbR^d} \Big[ b_{k_2}\cd \pa_x \d_m U_1^{k_1,k_2} + {1\over 2} \si^2 : \pa_{xx}^2 \d_m U^{k_1,k_2}_1 + f_{k_1} \Big] m(dx,1) =0,
\eea
and $U^{k_1,k_2}$ is Lipschitz continuous in $m$ under $\cW_1$ with a Lipschitz constant $C_{L, L_{k_1}}$ independent of $k_2$. Here, $\d_m U_1^{k_1,k_2}$ is in the sense of \reff{dmu1}. This, in particular, implies $|\pa_x \d_m U_1^{k_1,k_2}|\le C_{L, L_{k_1}}$ for all $k_2$. Then, we deduce from \eqref{eqdiff} that
\bea
\label{Uk2est}
\Big| (\dbL U^{k_1,k_2} + F_{k_1})(t, m) \Big| = \Big| \int (b - b_{k_2}) \cdot \pa_x \d_m U_1^{k_1,k_2} m(dx,1) \Big|   \le C_{L, L_{k_1}} \d_2,~\forall k_2 \ge 1,
\eea
where $F_{k_1}(t,m) := \int_{\dbR^d} f_{k_1}(t,x,m)m(dx,1)$ as in \reff{F}. Moreover, since
$$ U^{k_1, k_2}(t,m) =  g_{k_1}(\bar \dbP_{Y_T}^{t,m,k_2}) + \int_t^T F_{k_1}(r, \bar \dbP_{Y_r}^{t,m,k_2}) dr,$$
where $\bar \dbP^{t,m,k_2}$ is s.t. $X$ is unstopped with drift coefficient $b_{k_2}$ instead of $b$, one can easily show that
\beaa
\Big|U^{k_1,k_2}(t,m) -\big(g(\bar \dbP_{Y_T}^{t,m}) + \int_t^T F(r, \bar \dbP_{Y_r}^{t,m}) dr \big)\Big|\le C[\d_1 + \d_2] \le {\e\over 4n},
\eeaa
for $\d_1, \d_2$ small enough. Then
\beaa
V_n(t_{n-1}, m_{n-1}) = \sup_{m' \preceq m_{n-1}} \Big\{ \int_{t_{n-1}}^T F(r, \bar \dbP_{Y_r}^{t_{n-1},m'}) dr + g(\bar \dbP_{Y_T}^{t_{n-1},m'} ) \Big\} \le \sup_{m' \preceq m_{n-1}} U^{k_1, k_2}(t_{n-1}, m') +{\e\over 4n}.
\eeaa
By the supersolution property, $v$ is nondecreasing for $\preceq$, hence 
\beaa
\frac \e n \le (V_n  - v)(t_{n-1}, m_{n-1})  \le \sup_{m'\preceq m_{n-1}} (U_{k_1,k_2}-v)(t_{n-1}, m')  +{\e\over 4n}. 
\eeaa
This implies that, 
\bea\label{supersol-test}
\max_{(s, \tilde m)\in \cN_{\frac T n}(t_{n-1}, m_{n-1})} \Big\{(U^{k_1,k_2} - v)(s, \tilde m) - {T-s\over n}  \Big\} \ge \frac{3\e}{4n} - {T \over n^2}  \ge {\e \over 2n},  
\eea
for $n$ sufficiently large. Note that $(U^{k_1,k_2} - v)(T, \dbP_{Y_T}) \le (g_{k_1} - g)(\dbP_{Y_T}) \le {\e \over 4n}$ for all $\dbP \in \cP(t_{n-1}, m_{n-1})$ and $v$ is $\cN$-LSC, then by compactness of $\cN_{\frac T n}(t_{n-1}, m_{n-1})$ there exists an optimal argument $(t^*, m^*), \ t^* < T$, to the above maximum. Thus $\f(s, \tilde m) := U^{k_1,k_2}(s,\tilde m) - {T-s\over n} \in \ol \cA v(t^*,m^*)$, and therefore,
\beaa
0 &\le& -(\dbL \f + F)(t^*,m^*) = -(\dbL U^{k_1,k_2} + F_{k_1})(t^*,m^*)  + (F_{k_1} - F)(t^*,m^*) - {1\over n} \\
&\le& C_{L, L_{k_1}}\d_2 +  (F_{k_1} - F)(t^*,m^*)  - {1\over n},
\eeaa
where the last inequality thanks to \reff{Uk2est}.
Fix $k_1$  so that $(F_{k_1} - F)(t^*,m^*) \le \frac{1}{2n}$ and set $\d_2$ small enough, we obtain the desired contradiction. 
\ms

{\it Step 2:} We show that $(V_n - v)(t_{n-2}, \cdot) \le {2\e \over n}$. Assume to the contrary that there exists $m_{n-2}$ such that $(V_n-v)(t_{n-2}, m_{n-2}) > {2\e \over n}$.  By the DPP, we have
$$V_n(t_{n-2}, m_{n-2}) = \sup_{\dbP \in \cP_n(t_{n-2}, m_{n-2})} \Big\{ \int_{t_{n-2}}^{t_{n-1}} F(r, \dbP_{Y_r})dr + V_n(t_{n-1}, \dbP_{Y_{(t_{n-1})-}})\Big\},$$
Observe the fact that $v$ is a viscosity supersolution of \eqref{obstacle} also implies that $v + \frac \e n$ is a viscosity supersolution. Moreover, by {\it Step 1}, we have $(v + \frac \e n)(t_{n-1}, \cdot) \ge V_n(t_{n-1}, \cdot)$. Thus, using the same procedure as {\it Step 1} (where $V_n$ replaces $g$ on $(t_{n-2}, t_{n-1}]$), it follows that 
$$ \big(V_n - (v+\frac \e n)\big)(t_{n-2}, \cd) \le \frac \e n. $$ 
Finally, by backwards induction, we have $(V_n - v)(t_{n-j}, \cdot) \le {j\e \over n}$ for all $j\in \{0, \dots, n\}$, and thus
$ (V_n - v)(t, \cdot) \le \e, $
which implies by the arbitrariness of $n$ and $\e$ that $v \ge V$.
\qed

\subsection{Infinite horizon case}

As in \cite[\S 6.1]{TTZ}, we may formulate the problem in infinite horizon (i.e., in the case $T = \infty$), by replacing Assumption \ref{assum-bsig} with the following conditions:
\begin{assum}
\label{assum-infinity}
{\rm (i)} Assumption \ref{assum-bsig} holds true on $[0, \infty)$;\\
{\rm (ii)}  $\int_0^\infty \sup_{m\in \cP_2(\bS)} |F(t, m)|dt <\infty$;\\
{\rm (iii)}  For any $(t,m)$ and $\dbP\in \cP(t,m)$, $X_\infty:= \lim_{t\to \infty} X_t$ exists, $\dbP$-a.s.
\end{assum}
We remark that one sufficient condition of (ii) above is that $|f(t,x, m)|\le Ce^{-\l t}$ for some constants $C, \l>0$, and a special case of (iii) is:
\bea
\label{GBM}
d=1,\q b = b_0 x, \q \si = \si_0 x,\q b_0 - {1\over 2}\si_0^2 <0,
\eea
see e.g.\ Pedersen \& Peskir \cite{PP} and Xu \& Zhou \cite{XZ}. The last condition implies that, under $\bar \dbP^{t,m}$ in \reff{barP}, the unstopped process $X$ is a Geometric Brownian motion vanishing at infinity. 

Assumption \ref{assum-infinity} allows to include the case $\t = \infty$ in our framework and to preserve the compactness of $\cP(t,m)$ in the infinite horizon setting, so that all our previous results extend immediately.  

\begin{rem}{\rm
A study of the general infinite horizon would be of course of very relevant interest. In standard optimal stopping, this is addressed by adding a discount factor to the reward function $\dbE\big[e^{-r\t}\psi(X_\t)\big]$, assuming $F = 0$ for simplicity. However, embedding this in our formulation is more involved as $ \dbE\big[e^{-r\t}\psi(X_\t)\big] = \dbE\big[e^{-r\int_0^T I_s ds}\psi(X_\t)\big]$,
which is a function of the joint law of $X_\infty$ and the path of $I$. We therefore leave it for further research.}
\end{rem}

\section{Examples}
\label{sect-examples}

In this section we revisit the three examples studied in \cite{TTZ}, and add a new example concerning probability distortion. Note that in \cite{TTZ} we assumed that the value functions are smooth, which is hard to verify. In this section we show that they are the unique continuous viscosity solution of the corresponding obstacle problem. Note that we shall allow both $T<\infty$ and $T=\infty$, and correspondingly we always assume Assumption \ref{assum-bsig} or  \ref{assum-infinity}, and we shall report the detailed arguments in the case $T<\infty$ only. Moreover, for simplicity in this section we always assume $f=0$.

\subsection{Connection with standard optimal stopping}

Assume for this example that $b$ and $\si$ do not depend on the measure variable $m$. For a measurable function $\psi: \dbR^d\to\dbR$, we define the optimal stopping problem
\bea\label{MFstandard}
V(t,m) := \sup_{\dbP \in \cP(t,m)} \dbE^\dbP\Big[\psi(X_T)\Big],
&(t,m) \in \ol \bQ_0.&
\eea
That is, $g(\mu) := \int_{\dbR^d} \psi(x) \mu(dx)$ for $\mu\in \cP_2(\dbR^d)$. 
We  also introduce  $v(t,x):= V(t, \d_{(x,1)})$, which is related to  the standard obstacle problem: recalling \reff{cLx},
\bea\label{varineq}
\min\{-(\pa_t+\cL) v, v-\psi \} = 0, \q v(T,\cd) = \psi,\q\mbox{where}\q \cL v:= b\cd \pa_x v + {1\over 2} \si^2: \pa_{xx}^2 v.
\eea

\begin{prop}
\label{prop-standard}
 Assume $b, \si$ do not depend on $m$, $\si$ satisfies the regularities required in Lemma \ref{lem-Ureg}, and $\psi$ is uniformly continuous. Then $V$ is the unique continuous viscosity solution of the corresponding obstacle equation \reff{obstacle}, and it holds
\bea
\label{Vv}
V(t,m) =  \int_{\mathbf{S}}\big[v(t,x)i + \psi(x)(1-i)\big]m(dx,di).
\eea 
Moreover, there exists a pure strategy optimal stopping time.
\end{prop}
\proof First, by the uniform continuity of $\psi$ one can easily show that  $g$ is uniformly continuous in $m$ under $\cW_1$. Then by Theorem \ref{thm-reg} $V$ is continuous in $t$ and uniformly continuous in $m$ under $\cW_1$. Thus it follows from Theorems \ref{thm-existence} and \ref{thm-comparison} that $V$ is the unique viscosity solution of \reff{obstacle}. 

It remains to verify \reff{Vv}. Let $\dbP^*\in \cP(t, m)$ be such that, 
 \bea
 \label{standardtau}
 \t =  \inf\{s \ge t : v(s,X_s) =\psi(X_s) \},\q \dbP^*-\mbox{a.s. on} ~\{I_{t-}=1\}.
\eea
By the standard optimal stopping problem, see e.g. Karatzas \& Shreve \cite[Appendix D]{KS2}, $v$ is continuous and $\dbP^*$ is optimal. Then
 by \reff{MFstandard} we derive \reff{Vv}:
\beaa
V(t, m) =  \int_{\dbR^d} \!\!\!  \psi(x) m(dx, 0) + \int_{\dbR^d}\!\!\!  \dbE^{\dbP^*}\Big[\psi(X_T)\big| X_t = x\Big] m(dx, 1) = \int_{\dbR^d}\!\!\!  \psi(x) m(dx, 0) \!+\!\!  \int_{\dbR^d}\!\!\! v(t,x) m(dx, 1) . 
 \eeaa
Moreover, clearly the optimal stopping time determined by \reff{standardtau} is a pure strategy.
\qed

We remark that, by utilizing \reff{Vv}, it is possible to prove the uniqueness of viscosity solution under weaker requirement on $\si$. We leave the details to the interested readers.

\subsection{A generalization of the mean variance problem}
\label{cvx-example}
Consider the optimal stopping problem:
\bea\label{convex}
V(t, m) := \sup_{\dbP \in \cP(t,m)} \f\Big(\dbE^\dbP[\psi(X_T)]\Big), 
\eea
where $\psi: \dbR^d \to \dbR^k$ for some $k\ge 1$ and $\f: \dbR^k\to \dbR$. That is, $g(\mu) = \f\Big(\int_{\dbR^d} \psi(x) \mu(dx)\Big)$.

\begin{prop}
\label{prop-convex}
Let $b, \si$ satisfy the conditions in Theorem \ref{thm-comparison} {\rm (ii)}, $\psi$ be uniformly continuous and $\f$ be  continuous. Assume further that either $\psi$ is bounded or $\f$ is uniformly continuous.  Then $V$ is the unique continuous viscosity solution of the corresponding obstacle equation \reff{obstacle}.
\end{prop}
\proof Note that, when $|\psi|\le C$, we have $\big|\dbE^\dbP[\psi(X_T)]\big|\le C$ and thus in \reff{convex} we may replace $\f$ with the truncated function $\f_C(z):= \f({C\over |z|\vee C} z)$, $z\in\dbR^k$, which is uniformly continuous. Then in both cases, we may assume w.l.o.g. that $\f$ is uniformly continuous, and therefore, $g$ is uniformly continuous in $\mu$ under $\cW_1$. Then the results follows from Theorems \ref{thm-reg}, \ref{thm-existence}, and \ref{thm-comparison}.
\qed

\begin{rem}
\label{rem-convex} {\rm
 (i) In the case that $\f$ is convex: $\f(z) := \sup_{\a}[ \a z - \f^*(\a)]$, we have
\bea
\label{VVa}
V(t,m) = \sup_\a \ [V_\a(t,m) - \f^*(\a)],\q \mbox{where}\q V_\a(t,m) := \sup_{\dbP\in \cP(t, m)} \dbE^\dbP[ \a\cd \psi(X_T)]. \nonumber
\eea
Let $\a^*(t, m)$ be the optimal argument, then the optimal $\dbP^*$ for $V_{\a^*(t,m)}(t, m)$ is also  optimal for $V(t, m)$, and thus by Proposition \ref{prop-standard} there exists a pure optimal strategy for $V(t,m)$.

Moreover, let $\dbP^*$ be the optimal control for $V(0, m)$ and $V_{\a^*(0, m)}(0, m)$ as above, and denote $m^*_t:= \dbP^*_{Y_t}$. Then, by the DPP for $V$ and for $V_{\a^*(0,m)}$, we have
\bea
\label{a*DPP}
V(t, m^*_{t-}) = V(0, m) = V_{\a^*(0,m)}(0,m) - \f^*(\a^*(0,m)) = V_{\a^*(0,m)}(t,m^*_{t-}) - \f^*(\a^*(0,m)). \nonumber
\eea
That is, $\a^*(0,m)$ is optimal for $\sup_\a \ [V_\a(t,m^*_{t-}) - \f^*(\a)]$, or say, $\a^*(t, m^*_{t-}) = \a^*(0,m)$ for all $t$.

\no(ii) A more special case is the mean variance problem: for some constant $\l>0$,
\bea
\label{MV}
d=1, \q k=2, \q \psi_1(x) = x,\q \psi_2(x) = x^2,\q \f(z_1, z_2) = z_1 + {\l\over 2} z_1^2 - {\l\over 2} z_2.
\eea
In the homogeneous case \reff{GBM} with $T=\infty$, Pedersen $\&$ Peskir \cite{PP}  solved the problem $V(\d_{(x, 1)})$ and the optimal stopping time is a pure strategy. We are in a much more general framework. However, we should point out that \reff{MV} does not satisfies the technical conditions in Proposition \ref{prop-convex}.
\qed
}
\end{rem}

\subsection{Expected shortfall}

Let $d=1$, and fix some $\a \in (0,1)$, we consider the mean field optimal stopping problem
\bea\label{ESoptstop}
V(t,m) := \inf_{\dbP \in \cP(t,m)} \mathrm{ES}_\a^\dbP(X_T) \q \mbox{for all}
~~(t,m)\in\overline{\mathbf{Q}}_0, \nonumber
\eea
where $ \mathrm{ES}_{\a}^\dbP$ denotes the expected shortfall under $\dbP$, i.e., for any r.v. $Z$ with law $\mu$,
\bea\label{eq-RU}
&g(\mu):=\dis\mathrm{ES}_\a^\dbP(Z) := \frac 1 \a \int_0^\a q_\g(Z) d\g  =  \inf_{\b \in \dbR} \Big\{ \b + \frac{1}{1-\a}\int_\dbR (x-\b)^+\mu(dx)\Big\}, \\
&\dis\mbox{where}\q q_\g (Z) := \inf\{z : \mu(Z \le z) > \g \}. \nonumber
\eea
 Here the second equality has been established by Rockafellar \& Uryasev \cite{RU}.
 
 \begin{prop}
 $V$ is the unique continuous viscosity solution of the corresponding equation.
 \end{prop}
 \proof
 Clearly, $x \mapsto (x-\b)^+$ is Lipschitz continuous with Lipschitz constant $1$. By \eqref{eq-RU}, this implies that $g$ is Lipschitz continuous, and, given our assumptions on the coefficients, we conclude similarly to Proposition \ref{prop-convex} that the required claim follows. 
 \qed
 
Note further that
 \beaa
 V(t,m) = \inf_{\b\in \dbR} \Big\{ \b + \frac{1}{1-\a}V_\b(t,m)\Big\},\q \mbox{where}\q V_\b(t,m) :=  \inf_{\dbP\in \cP(t,m)} \dbE^\dbP[(X_T-\b)^+].
\eeaa
 One can easily show that $\lim_{\b\to\infty} \Big[\b + \frac{1}{1-\a}V_\b(t,m)\Big] = \lim_{\b\to -\infty} \Big[\b + \frac{1}{1-\a}V_\b(t,m)\Big] = \infty$, where the second equality is due to $\a\in (0,1)$. Then there exists optimal $\b^*=\b^*(t,m)\in \dbR$ such that  $V(t,m) =\b^* + \frac{1}{1-\a}V_{\b^*}(t,m)$. Therefore, similar to Remark \ref{rem-convex} (i), $V(t,m)$ and $V_{\b^*}(t,m)$ share an optimal $\dbP^*\in \cP(t, m)$, which is a pure optimal strategy as in Proposition \ref{prop-standard}.
 
 Moreover, in the homogeneous case with \reff{GBM} and $T=\infty$, one can easily show that $V$ and $V_\b$ are independent of $t$, and $V_\b(m)=\b^- m(\dbR_+,1) + \int_0^\infty (x-\b)^+ m(dx, 0)$ whenever $m(\dbR_+, \{0,1\}) =1$.

\subsection{Probability distortion}
Consider the following optimal stopping problem under probability distortion:
\bea
\label{distortion}
V(t, m):=\sup_{\dbP \in \cP(t,m)} \int_0^\infty \f\big(\dbP( \psi(X_T) \ge z)\big)dz,
\eea
where $\psi: \dbR^d \to [0, \infty)$ is a utility function, and $\f: [0,1]\to [0,1]$ is a probability distortion function: $\f(0)=0, \f(1)=1$, and $\f$ is strictly increasing. That is,
$
 g(\mu) =  \int_0^\infty \f\big(\mu(\{\psi\ge z\})\big)dz.
$

\begin{prop}
\label{prop-distortion}
Let $b, \si$ satisfy the conditions in Theorem \ref{thm-comparison} {\rm (ii)}, $\f$ be a uniformly Lipschitz continuous probability distortion function, and $\psi$ be uniformly continuous. Then $V$ is the unique continuous viscosity solution of the corresponding obstacle equation \reff{obstacle}.
\end{prop}
\proof As in the previous examples, it suffices to show that $g$ is uniformly continuous in $m$ under $\cW_1$.  Given arbitrary $\mu_1, \mu_2\in \cP_2(\dbR^d)$ and, for $i=1,2$, let $\xi_i$ be a random variable on $(\O, \cF, \dbP)$ such that $\dbP_{\xi_i} = \mu_i$ and $\dbE^{\dbP}[|\xi_1-\xi_2|]=\cW_1(\mu_1, \mu_2)$. Then
\beaa
\big|g(\mu_1) - g(\mu_2) \big| &\le& \int_0^\infty \Big| \f\big(\dbP( \psi(\xi_1) \ge z)\big) - \f\big(\dbP( \psi(\xi_2) \ge z)\big)\Big| dz \\
&\le& C\int_0^\infty \dbE^\dbP\Big[ \big| \1_{\{\psi(\xi_1) \ge z\}} -  \1_{\{\psi(\xi_2) \ge z\}}\big| \Big] dz
=  C\dbE^\dbP\Big[ \big|\psi(\xi_1) - \psi(\xi_2)\big|\Big].
\eeaa
Since $\psi$ is uniformly continuous, we see that $g$ is uniformly continuous in $\mu$ under $\cW_1$.
\qed

\begin{rem}
\label{rem-distortion2} {\rm
(i) In the homogeneous case \reff{GBM}  with $T=\infty$, Xu $\&$ Zhou \cite{XZ} solved the optimal stopping problem $V(\d_{(x,1)})$ for appropriate $\f, \psi$, and the optimal stopping time is a pure strategy.

\no(ii) The mean variance and probability distortion problems are typically viewed as time inconsistent, as the DPP does not hold for value function $v(t,x) := V(t, \d_{(x,1)})$.  However, we emphasize that, by viewing $m$ as our variable, $V$ satisfies the DPP and the problem is hence time consistent.
\qed}
\end{rem}

\appendix

\section{Technical results}\label{appendix}

 {\bf Proof of Lemma \ref{lem:compact-smallest}.} (i) The set $\{m':m' \preceq m\}$ is in continuous bijection with the compact set $\{ \hat m \in \cP_2(\mathbf{S} \times \{0,1\}) : \hat m \circ (\mathbf{x}, \mathbf{i})^{-1} = m \}$, with $(\mathbf{x}, \mathbf{i}, \mathbf{i'})$ the projection coordinates on $\mathbf{S} \times \{0,1\}$. This shows the compactness of $\{m':m' \preceq m\}$.
\\ 
(ii) As the map $m' \in \cK(m) \longmapsto m' (\dbR^d,1)$ is continuous and $\cK(m)$ is compact, there exists $\bar m \in \cK(m)$ s.t. $\bar m \in \underset{m' \in \cK(m)}{\mathrm{argmin}} m' (\dbR^d,1)$. 
Let $m' \in \cK(m)$ be such that $m' \preceq \bar m$ with some corresponding transition probability $p$, see the definition in \eqref{order}. Then, clearly $m' (\dbR^d,1) \le \bar m(\dbR^d,1)$ and thus equality holds by minimality of $\bar m(\dbR^d,1)$. As $p\le 1$, we conclude that $m' = \bar m$. 
\qed

\bs

\no {\bf Proof of Lemma \ref{lem-localbound}.}  For each $(t, m)\in \ol \bQ_0$ and $\dbP\in \cP(t, m)$, we extend $\dbP$ to $(\O, \cF_T)$ as follows:  denoted as $\hat\dbP\in \hat\cP(t, m)$,
$X_s = X_r$, $I_s = I_{t-}$, $s\in [-1, t)$, $\hat\dbP$-a.s.
We prove the lemma in two steps.

{\it Step 1.} For any compact $\cM\subset \cP_2(\bS)$, denote $\hat\cP_\cM:= \bigcup_{(t, m)\in [0, T]\times \cM} \hat\cP(t, m)$. For each $(t, m)$, $\dbP\in \cP(t, m)$, and $R>1$, following the proof of \cite[Proposition 2.2]{TTZ} we have
\beaa
&&\dis \dbE^{\hat\dbP}\big[  |X_T^*|^2 \big] \le C\int_{\dbR^d} |x|^2 m(dx,\{0,1\});\\
&&\dis \dbE^{\hat\dbP}\big[ |X^*_T|^2\1_{\{|X^*_T| \ge R\}} \big] \le C\int_\bS [1+|x|^2] \Big[\1_{\{|x|\ge \sqrt{R}-1\}}  +{1\over \sqrt{R}} \Big]m(dx, \{0,1\}). 
\eeaa
where $X^*_T:= \sup_{0\le s\le T} |X_s|$.
By the compactness of $\cM$, one can easily see that
\beaa
\sup_{\hat\dbP\in \hat\cP_\cM} \dbE^{\hat\dbP}\big[  |X_T^*|^2 \big]  <\infty,\q \lim_{R\to\infty} \sup_{\hat\dbP\in \hat\cP_\cM}\dbE^{\hat\dbP}\big[ |X^*_T|^2\1_{\{|X^*_T| \ge R\}} \big] =0.
\eeaa
Then the set $\hat\cP_\cM$ is compact. That is, for any $(t_n, m_n)\in [0, T]\times \cM$ and $\dbP^n\in \cP(t_n, m_n)$, there exists a subsequence,  still denoted the same, such that $\hat\dbP^n\to \hat \dbP^*$ under $\cW_2$, for some $\dbP^* \in \cP_2(\O, \cF_T)$. 

We may assume w.l.o.g. that $(t_n, m_n) \to (t^*, m^*)$ under $\cW_2$ for some $(t^*,m^*) \in   [0, T]\times \cM$. We next show that $\hat \dbP^* \in \hat\cP(t^*, m^*)$. Indeed, for any $\d>0$, we have $t^*-\d < t_n < t^*+\d$ for all $n$ large enough. By the required convergence, it is obvious that $X_s = X_{t^*-\d}, I_s = I_{t^*-\d}, s\le t^*-\d, \hat \dbP^*\mbox{-a.s.}$, and $\hat \dbP^*_{Y_{t-\d}} = m^*$. Thus, by sending $\d\to 0$, $X_s = X_{t^*}$, $I_s = I_{t^*-}$, $s < t^*$, $\hat \dbP^*$-a.s. and $\hat \dbP^*_{Y_{t^*-}} = m^*$. 
Here we used the fact that $X$ has continuous paths. Moreover,  following the arguments in \cite[Proposition 2.2]{TTZ} again, we see that the processes $M$ and $MM^\top$ in \reff{martingalepb} are $\hat\dbP^*$-martingales, on $[t^*+\d, T]$ for all $\d>0$, and hence also on $[t^*, T]$ (again since $X$ is continuous). That is, $\hat \dbP^* \in \hat\cP(t^*, m^*)$.

{\it Step 2.} We now show that $V$ is USC. Fix $(t, m)$ and choose $(t_n, m_n)\to (t, m)$ such that $\dis\lim_{n\to\infty} V(t_n, m_n) = \underset{(\tilde t, \tilde m)\to (t, m)}{\lim \sup} V(\tilde t, \tilde m)$. For each $n$, let $\dbP^n\in \cP(t_n, m_n)$ be optimal:
$
V(t_n, m_n) =  \int_{t_n}^T  F(r, \dbP^n_{Y_r})dr   + g(\dbP^n_{Y_T}).
$
Note that $\cM:= \{m, m_n, n\ge 1\}\subset \cP_2(\bS)$ is compact. By Step 1, we may assume without loss of generality that $\hat \dbP^n \to \hat \dbP \in \hat\cP(t, m)$. Then, since $F$ is continuous and $g$ is USC in $m$, we have
\beaa
\lim_{n\to\infty}V(t_n, m_n) = \lim_{n\to\infty}\Big[ \int_{t_n}^T  F(r, \hat\dbP^n_{Y_r})dr   + g(\hat\dbP^n_{Y_T})\Big] \le  \int_{t}^T  F(r, \hat\dbP_{Y_r})dr   + g(\hat\dbP_{Y_T})\le V(t,m).
\eeaa
 This means that $V$ is USC.
  \qed

\bs

\no {\bf Proof of Theorem \ref{thm-reg}.} (i) follows similar but easier arguments than (ii), so we prove (ii) only.   Let $\rho_0$ denote the modulus of continuity of $f, g$ under $\cW_1$. We proceed in two steps.

{\it Step 1.} Fix $t\in [0, T]$ and $m, \tilde m\in \cP_2(\bS)$. For any $\dbP\in \cP(t, m)$, by possibly enlarging the space, there exists $(\tilde X_t, \tilde I_{t-})$ on the space  $(\O, \cF, \dbP)$ such that 
\beaa
\dbP_{(\tilde X_t, \tilde I_{t-})} = \tilde m,~ \dbE^\dbP\Big[| \tilde X_t - X_t|+ |\tilde I_{t-} - I_{t-}|\Big] = \cW_1(m, \tilde m).
\eeaa
 Consider the following SDE on the space $(\O, \cF, \dbP)$:  for $\tilde Y := (\tilde X, \tilde I)$,
 \bea
\label{tildeX1}
\tilde X_s = \tilde X_t + \int_t^s b(r, \tilde X_r, \dbP_{\tilde Y_r})\tilde I_r dr + \int_t^s \sigma(r, \tilde X_r, \dbP_{\tilde Y_r)})\tilde I_r dW_r^\dbP, \q \tilde I_r:= I_r \tilde I_{t-},\q 
 \dbP\mbox{-a.s.}
\eea
Denote  $\D Y := \tilde Y-Y$.
Note that $I_r = I_r I_{t-}$, then 
\bea
\label{DI}
\sup_{t\le r\le T} |\D I_r| = I_r |\D I_{t-}|\le |\D I_{t-}|,\q\mbox{and thus}\q \dbE^\dbP\Big[\sup_{t\le r\le T} |\D I_r|\Big] \le  \cW_1(m, \tilde m).
\eea
Moreover,  for $\f=b, \si$, by the desired Lipschitz continuity under $\cW_1$, we have
\beaa
\Big|\f(r, \tilde X_r, \dbP_{\tilde Y_r)})\tilde I_r - \f(r, X_r, \dbP_{Y_r}) I_r\Big|  \le C\Big[|\D X_r| + \cW_1(\dbP_{\tilde Y_r}, \dbP_{Y_r})\Big] + C[1+|X_r|] |\D I_{t-}|.
\eeaa
By standard estimates, one can show that
\bea
&&\dis \dbE^\dbP_t\Big[\sup_{t\le s\le T} |X_s|^2\Big] \le C[1+|X_t|^2];\nonumber\\
&&\dis \dbE^\dbP_t\Big[|\D X_s|^2\Big] \le C\int_s^T \cW_1^2(\dbP_{\tilde Y_r}, \dbP_{Y_r})dr  + C|\D X_t|^2+C\sup_{t\le s\le T}\dbE^\dbP_t[1+|X_s|^2] |\D I_{t-}|^2;\nonumber\\
\label{DX}
&&\dbE^\dbP_t\big[ |\D X_s|\big] \le C\Big(\int_s^T \cW_1^2(\dbP_{\tilde Y_r}, \dbP_{Y_r})dr\Big)^{1\over 2} +C|\D X_t|+ C[1+|X_t|] |\D I_{t-}|.
\eea
This implies that
\beaa
&&\dbE^\dbP\big[ |\D X_s|\big] \le C\Big(\int_s^T \cW_1^2(\dbP_{\tilde Y_r}, \dbP_{Y_r})dr\Big)^{1\over 2} + C\dbE^\dbP\Big[|\D X_t|+[1+|X_t|] |\D I_{t-}|\Big];\\
 && \cW_1^2(\dbP_{\tilde Y_s}, \dbP_{Y_s}) \le C\int_s^T \cW_1^2(\dbP_{\tilde Y_r}, \dbP_{Y_r})dr + C\Big(\dbE^\dbP\Big[|\D X_t|+[1+|X_t|] |\D I_{t-}| \Big]\Big)^2.
  \eeaa
 By Grownwall inequality we have, for any $R>0$,
 \bea
 \label{dR}
 \sup_{t\le s\le T} \cW_1(\dbP_{\tilde Y_s}, \dbP_{Y_s}) &\le&  C\dbE^\dbP\Big[|\D X_t|+[1+|X_t|] |\D I_{t-}|\Big] \nonumber\\
 &\le& CR \cW_1(m, \tilde m) + C\dbE^\dbP\Big[|X_t|\1_{\{|X_t|\ge R\}}\Big]=: \d_R.
 \eea
  
Notice that $\tilde \dbP := \dbP \circ (\tilde X, I)^{-1} \in \cP(t, \tilde m)$. Then 
\bea
\label{DV1}
&&\dis \int_t^T \hat F(r, \dbP_{Y_r}) dr + g(\dbP_{Y_T}) - V(t, \tilde m)\nonumber \\
&&\dis \le \int_t^T \dbE^\dbP\big[f(r, X_r, \dbP_{Y_r})I_r  - f(r, \tilde X_r, \dbP_{\tilde Y_r})\tilde I_r \big]dr + \big[g(\dbP_{Y_T})  - g(\dbP_{\tilde Y_T}) \big]\\
&&\dis \le  \rho_0\Big( \cW_1(\dbP_{\tilde Y_T}, \dbP_{Y_T})\Big) + \int_t^T \dbE^\dbP\Big[\rho_0\big(|\D X_r|)+ \rho_0\big(\cW_1(\dbP_{\tilde Y_r}, \dbP_{Y_r})\big) +|f(r, X_r, \dbP_{Y_r})| |\D I_r|\Big].\nonumber
\eea
The uniform regularity of $f$ implies that
\beaa
|f(r, X_r, \dbP_{Y_r})| \le |f(r, 0, \dbP_{Y_r})| + C|X_r| \le C_m[1+|X_r|],
\eeaa
where the constant $C_m$ may depend on $m$. Then, by \reff{DI}, \reff{DX}, and \reff{dR}, we have
\beaa
&&\dbE^\dbP\Big[|f(r, X_r, \dbP_{Y_r})| |\D I_r|\Big] \le C_m \dbE^\dbP\Big[[1+|X_r|] |\D I_{t-}|\Big]\\
&&\dis \le C_m \dbE^\dbP\Big[|\D I_{t-}| + \rho_0(\d_R) + |\D X_t|+ [1+|X_t|]|\D I_{t-}|\Big] \le C_m [\cW_1(m,\tilde m)+ \rho_0(\d_R)\big].
\eeaa
Plug this into \reff{DV1}, we have
\beaa
\int_t^T \hat F(r, \dbP_{Y_r}) dr + g(\dbP_{Y_T}) - V(t, \tilde m) \le C_m \big[\cW_1(m,\tilde m)+ \rho_0(\d_R)\big] + \int_t^T \dbE^\dbP\big[\rho_0(|\D X_r|)] dr.
\eeaa
Since $\dbP\in \cP(t, m)$ is arbitrary, for some appropriate modulus of continuity  $\rho$ we have
$
V(t,m) - V(t, \tilde m) \le  C_m\rho(\d_R).
$
Switching $m, \tilde m$, and noticing that we may still use $X_t$ in $\d_R$, we have
\bea
\label{Step1est}
|V(t,m) - V(t, \tilde m)| \le  C_m\rho(\d_R).
\eea
Fix $m$ and send $\tilde m\to m$ under $\cW_1$, we see that 
\beaa
\underset{\tilde m\to m} {\lim \sup} \ |V(t,m) - V(t, \tilde m)| \le C_m\rho\Big(C\dbE^\dbP[|X_t|\1_{\{|X_t|\ge R\}}]\Big)
\eeaa
for any $R>0$. Now send $R\to\infty$, we see that $\lim_{\tilde m\to m} V(t, \tilde m) = V(t,m)$.

{\it Step 2.} Let $t<\tilde t$ and $m\in \cP_2(\bS)$.   By DPP we have
\bea
\label{DPP}
\left.\ba{c}
\dis V(t, m) = \sup_{\dbP\in \cP(t, m)} \Big\{ \int_t^{\tilde t} \!\! F(r, \dbP_{Y_r}) dr + V(\tilde t, \dbP_{Y_{\tilde t-}}) \Big\} = \sup_{\dbP\in \cP(t, m)} \Big\{ \int_t^{\tilde t} \!\! F(r, \dbP_{Y_r}) dr + V(\tilde t, \dbP_{Y_{\tilde t}})\Big\},\\
\dis V(\tilde t, m) = \sup_{m'\preceq m} V(\tilde t, m').
\ea\right.
\eea

First, for any $\dbP\in \cP(t, m)$, note that $m' := \dbP\circ (X_t, I_{\tilde t-})^{-1} \preceq m$, then
\beaa
V(\tilde t, \dbP_{Y_{\tilde t-}})  - V(\tilde t, m)  \le V(\tilde t, \dbP_{(X_{\tilde t}, I_{\tilde t-})})  - V(\tilde t, \dbP_{(X_t, I_{\tilde t-})}) \le C_m\rho(\d_R), 
\eeaa
thanks to \reff{Step1est} and \reff{dR}, where, following similar arguments as in Step 1,
\bea
\label{dR2}
\d_R:= CR \dbE^\dbP[|X_{\tilde t} - X_t|] + C\dbE^\dbP\Big[|X_t|\1_{\{|X_t|\ge R\}}\Big] \le CR\dbE^\dbP[1+|X_t|]\sqrt{\tilde t-t} + C\dbE^\dbP\big[ |X_t|\1_{\{|X_t|\ge R\}}\big].
\eea
Since $\dbP\in \cP(t, m)$ is arbitrary, by \reff{DPP} we have
\beaa
V(t, m)  - V(\tilde t, m)  \le \sup_{\dbP\in \cP(t, m)} \int_t^{\tilde t}  F(r, \dbP_{Y_r}) dr + C_m\rho(\d_R) \le C_m\rho(\d_R).
\eeaa

Next, for $m'\preceq m$, choose $\dbP\in \cP(t, m)$ s.t. $I_s = I_{t-}$, $t\le s< \tilde t$, and $\dbP\circ (X_t, I_{\tilde t})^{-1} = m'$. Then
\beaa
 V(\tilde t, m') - V(t, m) \le V(\tilde t, \dbP_{ (X_t, I_{\tilde t})}) - V(\tilde t, \dbP_{ (X_{\tilde t}, I_{\tilde t})}) - \int_t^{\tilde t}  F(r, \dbP_{Y_r}) dr\le C_m \rho(\d_R).
 \eeaa
Since $m'\preceq m$ is arbitrary, by \reff{DPP} we have
\beaa
V(\tilde t, m)  - V(t, m)  \le  C_m\rho(\d_R),\q\mbox{and thus}\q  \big|V(t, m)  - V(\tilde t, m) \big| \le  C_m\rho(\d_R).
\eeaa
This, together with \reff{dR2}, implies the desired regularity immediately.
\qed

\bs
\no{\bf Proof of Lemma \ref{lem-Ureg}.} We shall apply the results in Buckdahn, Li, Peng \& Rainer \cite{BLPR}. For this purpose, we extend functions on $\cP_2(\bS)$ to $\cP_2(\dbR^d\times \dbR)$. Let $\phi: \dbR\to \dbR$ be a smooth function with bounded derivatives s.t. $0\le \phi\le 1$, $\phi(0) = 0$, $\phi(1) = 1$, and $ \Phi: \hat m\in \cP_2(\dbR^d\times \dbR) \mapsto m \in \cP_2(\bS)$, with
\beaa
m(A, 1) := \int_\dbR \phi(y) \hat m(A, dy),\q m(A, 0) := \int_\dbR [1-\phi(y)] \hat m(A, dy),\q \mbox{for all $A\in \cB(\dbR^d)$}.
\eeaa
Now for $\f = b, \si, f, g$, define $\hat \f(t, x, \hat m) := \f(t,x, \Phi(\hat m))$. $\hat\f$ inherits the regularity of $\f$ on $\cP_2(\dbR^{d}\times \dbR)$. 

Next, fix a filtered probability space $(\hat \O, \hat \cF_T, \hat \dbF, \hat \dbP)$ on which is defined a $d$-dimensional Brownian motion $W$. For any $(t, \hat m)$, let $\xi\in \dbL^2(\cF_t; \dbR^d)$, $\eta\in \dbL^2(\cF_t; \dbR)$ be such that $\hat\dbP_{(\xi, \eta)} = \hat m$. Consider the following SDE on $[t, T]$ with solution $\hat Y = (\hat X, \hat I)$:
\beaa
\hat X_s = \xi + \int_t^s \hat b(r, \hat X_r, \hat\dbP_{\hat Y_r}) \phi(\hat I_r) dr +  \int_t^s \hat\si(r, \hat X_r, \hat\dbP_{\hat Y_r}) \phi(\hat I_r) dW_r;\q \hat I_s = \eta \1_{[t, T)}(s),\q \hat \dbP\mbox{-a.s.}
\eeaa 
We then define, recalling \reff{F},
\beaa
\hat U(t, \hat m) := \hat g(\hat\dbP_{\hat Y_T}) + \int_t^T \hat F(r, \hat\dbP_{Y_r})dr,\q \mbox{where}\q \hat F(r, \hat m) := \int_{\dbR^{d+1}} \hat f(r, x, \hat m) i ~ \hat m(dx, di).
\eeaa
We remark that, since $b$ and $\si$ are not necessarily bounded, the coefficients of the SDE for $\hat X$ is not Lipschitz continuous in $\hat I$. However, since $\hat I$ is already given, such Lipschitz continuity is not needed. In particular, we can apply \cite[Lemmas 6.2 and 7.1]{BLPR} so that 
$\pa_t \hat U, \pa_{\hat m} \hat U, \pa_{\hat y \hat m} \hat U$ exist and are continuous and bounded. Here $\pa_{\hat m} \hat U$ is the Lions derivative and satisfies:
$
\pa_{\hat m} \hat U(t, \hat m, \hat y) := \pa_{\hat y} \d_{\hat m} \hat U(t,\hat m, \hat y), 
$
see e.g. Carmona \& Delarue \cite[Vol. 1, Chapter 5]{CarDel}. We also remark that in \cite{BLPR} the function $\hat U$ takes the form $\hat U(t,x, \hat m)$ while here $\hat U$ does not have the $x$-variable. Moreover, note that each $m\in \cP(\bS)$ can be viewed as an element of $\cP(\dbR^d\times \dbR)$ with support included in $\bS$. Since $\phi(0) =0, \phi(1)=1$, one can easily see that $U(t, m) = \hat U(t, m)$. Then clearly $U\in C^{1,2}_2(\bQ_0)$. Finally, the $\cW_1$ Lipschitz continuity of $U$ follows arguments similar to that of Theorem \ref{thm-reg}, we thus omit the proof here.
\qed

\bs
\no{\bf Proof of Lemma \ref{lem-mol}} (ii) and (iii) follow directly from \cite[Theorem 3.1]{MZ}, after the straightforward extension to $\cP_2(\bS)$, as we will do next. Thus we shall only prove (i). For the ease of presentation, we assume $d=1$.

Fix $n \ge 1$, we construct $U_n$ as follows.   First, let $H_n, \phi^n_j\in C^\infty(\dbR)$, $j\in \dbZ$,  satisfy:
\beaa
& 0\le H_n \le 1,\q \supp(H_n) \subset [-{3n\over 2}, {3n\over 2}],\q H_n(x) = 1~\mbox{for}~ |x|\le n,\q |\pa_x H_n|\le {3\over n};\\
& 0\le \phi^n_j\le 1,\q \supp(\phi^n_j) \subset [{j-1\over n}, {j+1\over n}],\q  \phi^n_j (x) + \phi^n_{j+1}(x) =1 ~\mbox{for all}~ x\in [{j\over n}, {j+1\over n}].
\eeaa
See \cite[(3.3)]{WZ} for a construction of $\phi^n_j$.  Next,  for each $j\in \dbZ$, define
\bea
\label{psij}
\psi^n_j(\mu) := \int_\dbR \phi^n_j(x) H_n(x) \mu(dx) + \1_{\{j=0\}}  \int_\dbR [1- H_n(x)] \mu(dx),
\eea
for all finite measure $\mu$ on $\dbR$. We emphasize that, slightly different from \cite{WZ}, here the $\mu$ will be $m(\cd, i)$ whose total measure is less than $1$ and thus it is not a probability measure.  
Note that $\psi^n_j \ge 0$ and $\sum_{j\in \dbZ} \psi^n_j = \mu(\dbR)$. Moreover, denote $\dbZ_n := \{j\in \dbZ: |j|\le 2n^2\}$ with size $N_n := 4n^2+1$, and 
\beaa
\D_n := \Big\{ \vec z = \{z_j\}_{j\in \dbZ_n}: |z_j|\le N_n^{-3}~\mbox{for all}~j\neq 0,\q\mbox{and}\q z_0 := - \sum_{j\in \dbZ_n\backslash\{0\}} z_j\Big\}.
\eeaa
We now define, for each $\vec z\in \D_n$ and $m\in \cP_2(\bS)$, $i=0,1$,
\bea
\label{mny}
m_n(dx, i, \vec z) := \sum_{j\in \dbZ_n} \hat\psi_{j}^n(m(\cd, i), \vec z)\d_{ j\over n}(dx), \q\hat\psi^n_{ j}(\mu, \vec z) := {N_n \over N_n+1}\Big[ \psi^n_{ j}(\mu) + \mu(\dbR) [{1\over N_n^2} + z_{j}]\Big].
\eea
Note that $|z_0|\le N_n^{-2}$, and thus $\hat\psi^n_{ j}(\mu, \vec z)\ge 0$.  One may easily verify that 
\beaa
&\dis\sum_{j\in \dbZ_n} \hat\psi^n_{ j}(\mu, \vec z) = {N_n \over N_n+1}\Big[ \sum_{j\in \dbZ_n}\psi^n_{ j}(\mu) +  {\mu(\dbR)\over N_n} \Big] ={N_n \over N_n+1}\Big[ \mu(\dbR) +  {\mu(\dbR)\over N_n} \Big]=\mu(\dbR);\\
&\dis m_n(\dbR, i, \vec z)= \sum_{j\in \dbZ_n} \hat\psi_{j}^n(m(\cd, i), \vec z) = m(\dbR, i),\q\mbox{and thus}\q m_n(\bS, \vec z)=1.
\eeaa
In particular, this implies that $m_n(\cd, \vec z) \in \cP_2(\bS)$ for every $\vec z\in \D_n$, where the square integrability follows from the fact that $\supp(m_n(\cd, \vec z))$ is finite. 
Finally, let $\z_n$ be a smooth density function with support $\D_n$, and we construct
\bea
\label{Un}
U_n(m) := \int_{\D_n} U(m_n(\cd, \vec z)) \z_n(\vec z) d\vec z,\q m\in \cP_2(\bS).
\eea
The smoothness of $U_n$ follows from the same arguments as in \cite[Theorem 3.1]{MZ}. However, we note that \cite{MZ} uses the $\cW_1$-distance and requires $\cM$ to be a compact subset of $\cP_1(\bS)$. This is mainly for the uniform Lipschitz continuity of $U_n$ which holds only under $\cW_1$. Here we  provide a proof for the uniform convergence of $U_n$ under $\cW_2$.  We first show that 
\bea
\label{cMbar}
\overline \cM := \Big\{m_n(\cd, \vec z): m\in \cM, n\ge 1, \vec z\in \D_n\Big\} \subset \cP_2(\bS) \q \mbox{is compact}.
\eea
Indeed, fix $R>0$. Denote $\dbZ^R_{n}:= \{j \in \dbZ_n:  |j|\ge nR\}$ for  $n>{R\over 2}$. Then
\beaa
\dis \int_{\{|x|> R\}}\!\!\!\!\!\!\!\! |x|^2 m_n(dx, i, \vec z) 
= \sum_{j \in  \dbZ^R_{n}} {j^2\over n^2}\hat \psi^n_{\vec j}(m(\cd, i), \vec z) 
=  \sum_{j \in  \dbZ^R_{n}} {j^2\over n^2} {N_n \over N_n+1}\Big[ \psi^n_{j}(m(\cd, i)) + m(\dbR, i) [{1\over N_n^2} + z_{j}]\Big].
\eeaa
From the construction of $\psi_{\vec j}$, one can easily verify that
\beaa
\sum_{\vec j \in  \dbZ^R_{n}} {j^2\over n^2}  \psi^n_{j}(m(\cd, i)) \le 2 \int_{\{|x|> R\}} |x|^2 m(dx, i).
\eeaa
Moreover, note that $|z_{j}| \le N_n^{-3}$ for all $j\in \dbZ^R_{n}$. Then, for $n>{R\over 2}$,
\beaa
&& \int_{\{|x|> R\}} |x|^2 m_n(dx, i, \vec z) \le  2 \int_{\{|x|> R\}} |x|^2 m(dx, i) + \sum_{j \in  \dbZ^R_{n}} {j^2\over n^2}  m(\dbR, i) {C\over N_n^2}\\
 &&\le   2 \int_{\{|x|> R\}} |x|^2 m(dx, i) + {Cm(\dbR, i) \over N_n} \le 2 \int_{\{|x|> R\}} |x|^2 m(dx, i) +  { Cm(\dbR^d, i)\over R^2}.
\eeaa
On the other hand, when $n<{R\over 2}$, we have  $\int_{\{|x|> R\}} |x|^2 m_n(dx, i, \vec z)=0$. Thus, 
\beaa
\sup_{m\in \cM, n\ge 1, \vec z\in \D_n} \sum_{i=0,1} \int_{\{|x|> R\}} |x|^2 m_n(dx, i, \vec z) \le  2 \sup_{m\in \cM}\sum_{i=0,1}\int_{\{|x|> R\}} |x|^2 m(dx, i) +  { C\over R^2}.
\eeaa
Since $\cM\subset \cP_2(\bS)$ is compact, we have $\dis\lim_{R\to\infty}\sup_{m\in \cM}\sum_{i=0,1}\int_{\{|x|> R\}} |x|^2 m(dx, i)=0$. Then 
\beaa
\lim_{R\to\infty}\sup_{m\in \cM, n\ge 1, \vec z\in \D_n} \sum_{i=0,1} \int_{\{|x|> R\}} |x|^2 m_n(dx, i, \vec z) =0.
\eeaa
This   proves that $\overline \cM$ is uniformly square integrable, and therefore compact in $\cP_2(\bS)$.

Next, note that $\cM$ is also compact in $\cP_1(\bS)$, by \cite[(3.15)]{MZ} we have 
\bea
\label{W1conv}
 \lim_{n\to\infty} \sup_{m\in \cM, \vec z\in \D_n} \cW_1(m_n(\cd, \vec z), m) =0.
 \eea
  Then, for any $R>0$,  
\beaa
\cW_2^2(m_n(\cd, \vec z), m) \le R\cW_1(m_n(\cd, \vec z), m) + C\sum_{i=0,1} \int_{|x|\ge {R\over 2}} |x|^2\big[m_n(dx,i, \vec z)+m(dx,i)\big].
\eeaa
This, together with  the uniform integrability of $\overline \cM$ and \reff{W1conv}, implies immediately that
\bea
\label{W2conv}
 \lim_{n\to\infty} \sup_{m\in \cM, \vec z\in \D_n} \cW_2(m_n(\cd, \vec z), m) =0.
 \eea
Finally, by the compactness \reff{cMbar}, we see that $U$ is uniformly continuous on $\overline \cM$. Then it follows from \reff{Un} and \reff{W2conv} that  $\dis\lim_{n\to\infty} \sup_{m\in \cM}|U_n(m) - U(m)| =0$. \qed
 
\begin{rem}
\label{rem-Unmon}
{\rm 
While not used in the paper, the following property is interesting in its own right: if $U$ is monotone under $\preceq$, then so is the $U_n$ constructed in \reff{Un}. Indeed, assume $U$ is increasing, and let $m'\preceq m$ with transition probability $p$. For each $ \vec z\in \D_n$, by \reff{psij} and \reff{mny}, it is clear that
\beaa
0< \hat\psi^n_{ j}(m'(\cd, 1), \vec z) \le \hat\psi^n_{j}(m(\cd, 1), \vec z),\q\mbox{and thus}\q \hat p_{ j}(\vec z):= \hat\psi^n_{\vec j}(m'(\cd, 1),\vec z) ~\slash~ \hat\psi^n_{ j}(m'(\cd, 1),\vec z) \in (0,1].
\eeaa
Since $m'(dx, \{0,1\}) = m(dx, \{0,1\})$, it is also obvious that  
$
\dis\sum_{i=0,1}\hat\psi^n_{j}(m'(\cd, i), \vec z) =\sum_{i=0,1} \hat\psi^n_{j}(m(\cd, i),\vec z).
$
Then $m'_n(\cd, \vec z) \le m_n(\cd, \vec z)$ for each $\vec z\in \D_n$, with transition probability $\hat p(\cd, \vec z)$ satisfying $\hat p({j\over n}, \vec z) = \hat p_{j}(\vec z)$ for all $j\in \dbZ_n$. Then, since $U$ is increasing, by \reff{Un} we see that $U_n(m') \le U_n(m)$. 
\qed}
\end{rem}

\no {\bf Proof of Lemma \ref{lem-lusin}.} (i) 
Let $m' \preceq m$ with transition probability $p$. As $m$ is a probability measure on $(\bS, \cB(\bS))$, it is a Radon measure. Then, by Lusin's theorem (see Folland \cite[7.10]{Folland}), we may find for all $k \ge 1$ a continuous $p_k : \dbR^d \longrightarrow [0,1]$ s.t. 
\beaa
m\big( \{ x : p(x) \neq p_k(x) \}, \{0,1\}\big) \le \frac 1 k.
\eeaa
Let $\{m_k'\}_{k \ge 1}$ be the measures obtained from $m$ by applying the transition probabilities $\{p_k\}_{k \ge 1}$, and $\phi$ a bounded and continuous function. Then
\beaa
\Big\lvert \int_{\dbR^d } \phi(x)p_k(x)m(dx,1) - \int_{\dbR^d }\phi(x) p(x) m(dx,1) \Big\rvert 
\le
\frac 2 k \lVert \phi \rVert_{\infty} ,
\eeaa
and thus $m'_k(dx,1)$ converges weakly to $m'(dx,1)$. We do similarly with $m'_k(dx,0)$, and thus $m'_k$ converges weakly to $m'$. As $\{m'_k\}_{k \ge 1}$ is uniformly integrable, we have $\lim_{k \to \infty} \cW_2(m_k',m') = 0$. As $v$ is nondecreasing for $\preceq$, we have $v(m) \ge v(m'_k)$ for all $k \ge 1$. Then, as $v$ is $\cN$-LSC, we have 
$v(m) \ge \underset{k \to \infty}{\lim \inf} \ v(m'_k) \ge v(m').$

\no(ii) As $(\dbD_I \f)_*$ is LSC,
there exists $ \d > 0$ s.t. $(\dbD_I \f)_* \ge 0$ on $[t, t+ \d] \times \cB_{\cW_2}(m,\d)$. Let $(s, m_0), (s, m_1)$ be in this neighborhood, s.t. $m_1 \preceq m_0$ with transition probability $p$. Then, we have
\bea\label{loc-incr}
\f(s, m_0) - \f(s, m_1) = \int_0^1 \int_{\dbR^d} D_I \f(t, \l m_0 + (1-\l) m_1, x) (1-p(x)) m(dx,1) d\l. \nonumber
\eea
By convexity of $\cB_{\cW_2}(m,  \d)$, we have $D_I \f(t, \l m_0 + (1-\l) m_1, \cdot) \ge 0$, hence the desired result. 
\qed

\bibliographystyle{plain}
\bibliography{references2}
\end{document}